\documentclass[lettersize,journal]{IEEEtran} \def\TwoColumn{1}

\usepackage{xcolor}
\usepackage{amsmath,amssymb,amsfonts,amsthm}
\usepackage{algorithmic}
\usepackage{graphicx}
\usepackage{dsfont}
\usepackage{pstricks}
\usepackage{afterpage}
\usepackage{comment}
\usepackage{cite}
\usepackage{url}
\usepackage{etoolbox}
\usepackage{hyperref}
\hypersetup{hidelinks=true}
\usepackage{textcomp}
\def\BibTeX{{\rm B\kern-.05em{\sc i\kern-.025em b}\kern-.08em
    T\kern-.1667em\lower.7ex\hbox{E}\kern-.125emX}}
\allowdisplaybreaks[1]
\newcommand{\diag}{\mbox{diag}}
\newcommand{\eqdef}{\buildrel \triangle \over =}

\ifdef{\TwoColumn}{}{ \usepackage{setspace}\doublespacing}

\makeatletter  
\def\@endtheorem{$\blacksquare$\endtrivlist\@endpefalse } % insert `\qed` macro
\makeatother

\newtheoremstyle{thmstyle}%                % Name
{}%                                     % Space above
{}%                                     % Space below
{}%                                     % Body font
{}%                                     % Indent amount
{\bfseries}%                            % Theorem head font
{.}%                                    % Punctuation after theorem head
{ }%                                    % Space after theorem head, ' ', or \newline
{}%                                     % Theorem head spec (can be left empty, meaning `normal')
\theoremstyle{thmstyle}
\newtheorem{definition}{Definition}
\newtheorem{lemma}{Lemma}
\newtheorem{thrm}{Theorem}
\usepackage[overload]{textcase}

\usepackage{textcomp}
\def\BibTeX{{\rm B\kern-.05em{\sc i\kern-.025em b}\kern-.08em
    T\kern-.1667em\lower.7ex\hbox{E}\kern-.125emX}}

\begin{document}
\title{Prescribed-Time Stability Properties of Interconnected Systems}
\author{Prashanth Krishnamurthy, Farshad Khorrami, and Anthony Tzes
  \thanks{The first and second authors are with the Control/Robotics Research Laboratory (CRRL),
    Dept. of ECE, NYU Tandon School of Engineering, Brooklyn, NY 11201, USA, and with NYUAD Center for Artificial Intelligence and Robotics (CAIR), Abu Dhabi, UAE. The third author is with Dept. of ECE, NYU Abu Dhabi (NYUAD), Abu Dhabi, UAE and with NYUAD Center for Artificial Intelligence and Robotics (CAIR), Abu Dhabi, UAE.
    (e-mails: \{prashanth.krishnamurthy, khorrami, anthony.tzes\}@nyu.edu).}
  \thanks{This work is supported in part by the New York University Abu Dhabi (NYUAD) Center for Artificial Intelligence and Robotics, funded by Tamkeen under the NYUAD Research Institute Award CG010.}
  }

\maketitle

\bstctlcite{IEEEexample:BSTcontrol}

\begin{abstract}                % Abstract of not more than 250 words.
  Achieving control objectives (e.g., stabilization or convergence of tracking error to zero, input-to-state stabilization) in ``prescribed time'' has attracted significant research interest in recent years. The key property of prescribed-time results unlike traditional ``asymptotic'' results is that the convergence or other control objectives are achieved within an arbitrary designer-specified time interval instead of asymptotically as time goes to infinity. In this paper, we consider cascade and feedback interconnections of prescribed-time input-to-state stable (ISS) systems and study conditions under which the overall states of such interconnected systems also converge to the origin in the prescribed time interval. We show that these conditions are intrinsically related to properties of the time-varying ``blow-up'' functions that are central to prescribed-time control designs. We also generalize the results to interconnections of an arbitrary number of systems. As an illustrative example, we consider an interconnection of two uncertain systems that are prescribed-time stabilized using two different control design methods and show that the two separate controllers can be put together to achieve prescribed-time stability of the interconnected system. 
 % interconnection terms between
 % the systems and show that a control design can be performed to achieve
 % prescribed-time stabilization of the interconnected system.
\end{abstract}


\section{Introduction}

Prescribed-time stabilization/regulation has been increasingly attracting
interest in the controls literature over recent years
\cite{prescribed_time2,JSL17,TYS18,BVAD18,prescribed_time3,SDM19,ZSW19,KKK19b,KKK19a,KKK20a,KK20b,Gomez20,ZS21,KKK21a,HNL22,KK23,ZZ23,YS23,LK23}.
Unlike the classical ``asymptotic'' control designs \cite{KKK95,Isi99} in which the control
objective (e.g., stabilization to the origin) is sought to be attained over the
infinite time horizon (i.e., as time $t\rightarrow\infty$), prescribed-time
control seeks to achieve the control objective within an {\em a priori} chosen
constant time $T$ that is free to be picked by the designer irrespective of the
initial conditions.
Apart from the benefit of being able to specify a desired
convergence time, prescribed-time control is valuable in applications where the
control task is inherently defined over a specific time interval (e.g., in
applications such as autonomous vehicle rendezvous and missile guidance).
Prescribed-time control is closely linked to the notions of ``finite-time'' and
``fixed-time'' control. While the terminal time $T$ in prescribed-time control
is constant and independent of the initial conditions, finite-time control
\cite{finite_time1,finite_time2,finite_time3,finite_time4,finite_time5,finite_time6,finite_time7,finite_time9,APA08}
seeks instead to just enforce the control objective within a finite time (but not
necessarily constant or selectable by the designer) and fixed-time control
\cite{Pol12,PEP16,ZPEP18,AGJSD19} enforces that the
terminal time is constant (independent of the initial conditions), but not
necessarily arbitrarily specifiable by the designer. 
A broad overview of the literature on prescribed-time control and its links to
finite-time and fixed-time control is provided in the recent survey paper
\cite{SYL23}. 

In the prescribed-time control designs in the literature, two main strategies
can be enumerated:
\begin{itemize}
\item State scaling (e.g.,
  \cite{prescribed_time2,prescribed_time3}): The system state
  is scaled by a time-dependent function that grows (``blows up'') to infinity
  as $t\rightarrow T$ and the control design seeks to keep this scaled state
  bounded, thereby implicitly forcing the actual system state to go to zero. For
  example, defining a ``blow-up'' function $\mu(t)$ such that
  $\mu(t)\rightarrow\infty$ as $t\rightarrow T$ and defining the scaled state as
  $\tilde x = \mu(t)x$ where $x$ is the original system state, the control
  design seeks to keep $\tilde x$ bounded and therefore make $x$ go to zero as $t\rightarrow T$.
\item Time scaling (e.g., \cite{KKK19a,KKK19b,KK20b}): A nonlinear
  time scale
  transformation is introduced as $\tau = a(t)$. Picking $a$ to be a function
  such that $a(0)=0$ and $\lim_{t\rightarrow T}a(t)=\infty$, this temporal
  transformation maps the finite prescribed time interval $t\in[0,T)$ to the
  infinite time interval $\tau\in[0,\infty)$. Then, in terms of the new time
  variable $\tau$, the original prescribed-time control objective reduces to an
  asymptotic objective. Hence, control designs that address asymptotic
  convergence (e.g., dynamic high gain based control designs \cite{KK04f,KK06,KK08,KK16b}) can be applied to enforce the control objective as
  $\tau\rightarrow \infty$, therefore implicitly achieving prescribed-time
  properties as $t\rightarrow T$ in terms of the original time variable $t$.
\end{itemize}
A common aspect of both these strategies is the presence of time-dependent
functions that go to $\infty$ as $t\rightarrow T$, i.e., the blow-up functions
in state scaling and the time scale transformation functions in time scaling.
Further links between these methods can be identified such as the fact that the
derivative of a time scale transformation function is naturally a blow-up
function (and a time scale transformation can be constructed from the integral of
a blow-up function).
Another common property of both these strategies is that the control gain goes to
$\infty$ as $t\rightarrow T$.
Indeed, as noted in \cite{prescribed_time2,prescribed_time3},
this property that control gains go to infinity is to be expected for any
approach for regulation in prescribed finite time (including, for example,
optimal control based designs with a terminal constraint and sliding mode based
controllers with time-varying gains). However, the control input itself stays
bounded (and goes to zero under appropriate conditions) even though the control
gains go to infinity since the control law involves products of the control
gains and system state variables and the state variables go to zero fast enough
relative to how fast the control gains grow. 

Prescribed-time control designs have been developed under several scenarios over
recent years including both state-feedback \cite{KKK19a} and output-feedback
contexts \cite{KKK19b},
incorporating adaptation to uncertain parameters \cite{KKK20b,HNL22} that could be
coupled with unmeasured state variables \cite{KKK21a}, systems with unknown
input gain and appended dynamics \cite{KK23}, systems with unknown time delays
\cite{KK20b}, nonlinear control for linear systems in canonical form
\cite{YS23}, and linear time-varying feedback for a class of nonlinear systems
\cite{ZS21}. Along with development of control designs for
stabilization/regulation, extensions of input-to-state stability (ISS) concepts
\cite{sontag_iss} have been addressed in the finite-time/fixed-time
\cite{HJF10_FTISS,LEPP20} and prescribed-time  
\cite{prescribed_time2,prescribed_time3,ZZD22_PTISS} contexts.

In this paper, we draw inspiration from classical studies of interconnected
ISS systems such as the small-gain theorem \cite{jiang_iss} and ask the question of under what conditions
interconnections of prescribed-time ISS (PT-ISS) systems would be prescribed-time stable.
This question turns out to have surprising intricacies especially when the
systems in the interconnection have different blow-up functions (which as noted
above are intrinsic to prescribed-time control designs). 
The analysis of interconnected prescribed-time ISS systems necessitates several
new notions introduced in this paper including polynomially bounded blow-up
functions, prescribed-time exponential convergence, and prescribed-time
exponentially convergent Lyapunov certificates. These new notions also provide
insights into aspects of the prescribed-time control designs such as the
boundedness of the control signal due to the rapidity of state convergence to
zero being faster than the rapidity of control gain explosion to infinity. We
present several properties of these notions and study cascade and feedback
interconnections of prescribed-time ISS systems both in the context of two
systems in the interconnection and an arbitrary number of systems. We develop
sufficient conditions under which interconnections of an arbitrary number of
nonlinear systems retain prescribed-time convergence properties. We also present
simulation examples of cascade and feedback interconnections to illustrate the
developed concepts.

This paper is organized as follows. The basic terminologies and definitions are summarized in Section~\ref{sec:preliminaries}.
Then, in Section~\ref{sec:exponential_convergence}, the notions of polynomially bounded blow-up functions and prescribed-time exponential convergence and the links between these notions are introduced along with several lemmas.
The main results of this paper on cascade and feedback interconnections of
PT-ISS systems are presented in Section~\ref{sec:interconnections}. Illustrative
examples of cascade and feedback interconnections of two nonlinear systems are
considered in Section~\ref{sec:examples}. It is shown along with simulation studies that control laws designed separately for the two systems using two different prescribed-time control design methods can be put together to yield prescribed-time stability of the interconnected system as an application of the results developed in Section~\ref{sec:interconnections}. Concluding remarks are summarized in Section~\ref{sec:conclusion}.

\section{Preliminaries}
\label{sec:preliminaries}

Throughout this paper, we use the notation $T$ to denote the prescribed time ($T > 0$), which is the finite designer-specified time at which control objectives (e.g., stabilization, input-to-state stabilization) are desired to be achieved. Also, we will consider the initial time to be 0, i.e., the prescribed time interval is $[0,T)$.

%% The definitions of the prescribed-time properties are analogous to prior literature \cite{}. However, we view them from the perspective of time scale transformations (Definition~\ref{defn:temporal_transformation}) similar to those used in our prior prescribed-time controller designs in \cite{} and note that the Lyapunov characterizations have intuitive interpretations from the viewpoint of time scale transformations.

To outline the prescribed-time convergence and ISS properties, we will consider systems of the following two forms:
\begin{itemize}
\item time-varying systems without exogenous input signal:
\begin{align}
  \dot x &= f(x,t)
           \label{eq:sys_without_exogenous}
\end{align}
with $x\in{\cal R}^n$ being the state, $t\in {\cal R}$ the time, and $f$ a continuous function of its arguments.
\item time-varying systems with exogenous input signal $d$:
\begin{align}
\dot x &= f(x,d,t)
         \label{eq:sys_with_exogenous}
\end{align}
with $x\in{\cal R}^n$ being the state, $d\in{\cal R}^{n_d}$ being a uniformly bounded exogenous input signal, $t\in {\cal R}$ the time, and $f$ a continuous function of its arguments.
\end{itemize}

To formalize the prescribed-time stability and convergence properties of the
systems above, the concept of a blow-up function \cite{LK23} will be very useful as
summarized in the definition below followed by the definitions of key
prescribed-time stability and convergence properties.

\begin{definition}
  Blow-Up Function: A function $\varphi:[0,T)\rightarrow [0,\infty)$ is said to be a blow-up
  function if 
  it is monotonically increasing over $[0,T)$, continuously differentiable,
  $\varphi(t)<\infty$ for all $t<T$, $\varphi(t)\geq \varphi_0$ for all $t\in[0,T)$ with $\varphi_0$ being a positive constant, $\lim_{t\rightarrow
    T}\varphi(t)=\infty$, and $\lim_{t\rightarrow T}\int_0^t \varphi(s)ds =
  \infty$. An analogous function is also referred to as a
  $T$-finite-time-escaping function in \cite{ZZD22_PTISS,ZZ23}.
\end{definition}

\begin{definition}
  \label{defn:ptc}
  Prescribed-Time Convergent (PT-C) System: The system
  \eqref{eq:sys_without_exogenous} is said to be PT-C if starting from any
  initial condition $x(0)=x_0$, the inequality
    $|x(t)|\leq \beta(|x_0|,\mu(t)-\mu(0))$
  holds with $\beta$ being a class ${\cal KL}$ function\footnote{A continuous
    function $\alpha:[0,a)\rightarrow[0,\infty)$ is said to be of class ${\cal K}$
    if it is strictly increasing and $\alpha(0)=0$.  It is said to be of
    class ${\cal K}_\infty$ if furthermore $a=\infty$ and $\alpha(r)\rightarrow \infty$
    as $r\rightarrow\infty$. A continuous function
    $\beta:[0,\infty)\rightarrow[0,\infty)$ is said to be of class ${\cal L}$ if it is monotonically decreasing and
    $\lim_{s\rightarrow\infty}\beta(s)=0$.
    A class ${\cal KL}$
    function is class ${\cal K}$ with respect to its first argument and
    class ${\cal L}$ with respect to its second argument.} and $\mu$ being a blow-up function. PT-C implies that
  $\lim_{t\rightarrow T}|x(t)|=0$ starting from any initial condition.
\end{definition}

\begin{definition}
  \label{defn:pt_iss}
  Prescribed-Time ISS (PT-ISS) System \cite{prescribed_time2,prescribed_time3,ZZD22_PTISS}: The system \eqref{eq:sys_with_exogenous} is said to be PT-ISS if starting from any initial condition $x(0)=x_0$ and with any exogenous input signal $d:[0,T)\rightarrow {\cal R}$, the inequality
  $|x(t)|\leq \beta(|x_0|,\mu(t)-\mu(0)) + \gamma(\sup_{s\in[0,t)}|d(s)|)$
  holds with $\beta$ being a class ${\cal KL}$ function, $\mu$ being a blow-up function, and $\gamma$ being a class ${\cal K}$ function. PT-ISS implies that the effect of the initial condition $x_0$ goes to 0 as $t\rightarrow T$.
\end{definition}

\begin{definition}
  Prescribed-Time ISS + Convergent (PT-ISS+C) System \cite{prescribed_time2,prescribed_time3,ZZD22_PTISS}: The system \eqref{eq:sys_with_exogenous} is said to be PT-ISS+C if starting from any initial condition $x(0)=x_0$ and with any exogenous input signal $d:[0,T)\rightarrow {\cal R}$, the inequality
  $|x(t)|\leq \beta(\beta_x(|x_0|)+ \gamma(\sup_{s\in[0,t)}|d(s)|),\mu(t)-\mu(0))$
  holds with $\beta$ being a class ${\cal KL}$ function, $\mu$ being a blow-up function, and $\beta_x$ and $\gamma$ being class ${\cal K}$ functions. PT-ISS+C implies that $x$ goes to 0 as $t\rightarrow T$ despite the exogenous input signal $d$ which does not need to go to 0 as $t\rightarrow T$.
\end{definition}

The concept of a blow-up function is closely related to the concept of a time
scale transformation as introduced in the definition below to map the finite
time interval $[0,T)$ to the infinite time interval $[0,\infty)$.

\begin{definition}
  \label{defn:temporal_transformation}
  Time Scale Transformation \cite{KKK19a,KKK20a}: A function $a:[0,T)\rightarrow[0,\infty)$ is said to be a time scale transformation (over the prescribed time interval given by $[0,T)$) if 
  $a(0)=0$, $\lim_{t\rightarrow T}a(t)=\infty$, $a(t) < \infty$ for all $t < T$, 
twice continuously differentiable, 
  and $a'(t)\eqdef \frac{da}{dt}$ is a blow-up function.
\end{definition}

In Section~\ref{sec:interconnections}, a class of matrices referred to in the linear algebra literature as Lyapunov diagonally stable matrices \cite{BBP78_diagstability,BH83_diagstability,Kra91_diagstability} and the estimation of ``weighted decay rates'' of the asymptotically stable linear systems $\dot x = Ax$ with $A$ being a Lyapunov diagonally stable matrix will be seen to be required in the formulation of conditions for prescribed-time stability of interconnected systems. For this purpose, we introduce the definitions of these concepts below.

\begin{definition}
  Lyapunov Diagonally Stable Matrix\footnote{While this terminology implies that
    a Lyapunov diagonally stable matrix has all eigenvalues in the open right
    half plane and is therefore opposite to the usage of the word ``stable'' in
    the controls literature, this terminology is adopted here for consistency
    with the linear algebra literature 
    \cite{BBP78_diagstability,BH83_diagstability,Kra91_diagstability}.} \cite{BBP78_diagstability,BH83_diagstability,Kra91_diagstability}: A matrix $A$ is Lyapunov diagonally stable if a diagonal matrix $P>0$ exists such that $PA+A^TP > 0$.
\end{definition}

It is known from \cite{BBP78_diagstability} that if an $n\times n$ matrix $A$ is a Lyapunov diagonally stable matrix, then the following property holds: a $1\times n$ vector $q>_e 0$ exists\footnote{The notation $m \leq_e n$ with $m$ and $n$ being vectors (or matrices) of the
    same dimensions indicates element-wise inequalities between corresponding elements of $m$ and $n$. The notation $m\leq_e 0$ indicates that all elements of $m$ are non-positive.  The
    notations $\geq_e$, $<_e$, and $>_e$ are defined analogously.} such that $qA >_e 0$.

% \begin{definition}
%   Decay rate: Given a strict Hurwitz matrix $A$ (i.e., all eigenvalues in the open left half plane) of dimension $n\times n$ such that $-A$ is a Lyapunov diagonally stable matrix, the decay rate $\delta(A)$ is defined as the solution of the following optimization problem (where $y\in {\cal R}$ and $P$ a diagonal $n\times n$ matrix)\footnote{If $P$ is a square matrix, the notation $P>0$ denotes that $P$ is a symmetric positive-definite matrix. The notations $P\geq 0$, $P<0$, etc., are defined analogously.}:
%   \begin{align}
%     \max y \mbox{ subject to } P > 0, PA+A^TP < -2y P, P \mbox{ diagonal}.
%   \end{align}
% \end{definition}
% With $-A$ being a Lyapunov diagonally stable matrix, it is evident that the above optimization problem is well-defined and can, for example, be numerically solved using the generalized eigenvalue (gevp) function in Matlab. We denote a $P$ corresponding to the optimal $y$ as $P_\delta(A)$. Note that $P_\delta(A)$ is non-unique. For our purpose, $P_\delta(A)$ can simply be interpreted as any choice of $P$ corresponding to the optimal $y$.

\begin{definition}
  Weighted Decay Rate: Given a strict Hurwitz matrix $A$ (i.e., all eigenvalues in the open left half plane) of dimension $n\times n$ such that $-A$ is a Lyapunov diagonally stable matrix, the weighted decay rate $\delta(A)$ is defined as the solution of the following optimization problem (where $y\in {\cal R}$ and $q$ a $1\times n$ vector):
    $\max y \mbox{ subject to } q >_e 0, qA <_e -y q$.
\end{definition}
With $-A$ being a Lyapunov diagonally stable matrix, it is evident that the
above optimization problem is well-defined and can, for example, be numerically
solved using the generalized eigenvalue (gevp) function in Matlab. Note that in terms of $q$, the optimization problem involves linear matrix inequalities (LMIs). The scalar $y$ appears in a bilinear combination as $yq$.
Note that with $-A$ being a Lyapunov diagonally stable matrix, we will have $\delta(A)>0$.
We denote a $q$ corresponding to the optimal $y$ as $q_\delta(A)$. Note that $q_\delta(A)$ is non-unique. For our purpose, $q_\delta(A)$ can simply be interpreted as any choice of $q$ corresponding to the optimal $y$.

The motivation for viewing $\delta(A)$ as a weighted decay rate arises from the following observation. Considering the system $\dot x= Ax$ and considering that the states $x$ are constrained to be non-negative, define the scalar combination of the state variables in $x$ as $z=qx$. Then, we obtain $\dot z=qAx \leq -\delta(A)qx=-\delta(A)z$. Hence, $z$ exponentially converges to 0 with an exponential rate given by $\delta(A)$, therefore motivating the terminology ``weighted decay rate.''

\section{Polynomial Boundedness and Prescribed-Time Exponential Convergence}
\label{sec:exponential_convergence}

A time scale transformation $\tau=a(t)$ as defined in Section~\ref{sec:preliminaries} maps the finite time interval $[0,T)$ in terms of the original time variable $t$ to the infinite time interval $[0,\infty)$ in terms of the transformed time variable $\tau$. Hence, prescribed-time control objectives formulated over the finite time interval $[0,T)$ in terms of the original time variable $t$ can be viewed as asymptotic control objectives formulated over the infinite time interval $[0,\infty)$ in terms of the transformed time variable $\tau$.
Since a time scale transformation is monotonically increasing by the properties in its definition, it is seen that such a function is invertible. We denote the inverse function of the time scale transformation $a(.)$ by $a^{-1}(.)$. The derivative $a'(t)$ can be written in terms of the transformed time variable $\tau$ as $\tilde a'(\tau)=a'(a^{-1}(\tau))$.
Note that if $a$ is a time scale transformation, then $a'(t)=\frac{da}{dt}$ is a
blow-up function.
Furthermore, if $a$ is a time scale transformation, then $a(t)+a_0$ is a blow-up
function with any positive constant $a_0$.
Also,  if $\varphi$ is a blow-up function, then it is seen that the function defined as $a_\varphi(t)=\int_0^t \varphi(s)ds$ is a time scale transformation.
An important subset of blow-up functions that is widely used in the literature \cite{prescribed_time2,prescribed_time3} is based on the expression $\frac{T}{T-t}$ that goes to infinity as $t\rightarrow T$. Based on this expression, a class of blow-up functions and time scale transformations is defined below.

\begin{definition}
  \label{defn:polynomially_bounded}
  Polynomially Bounded Blow-Up Function:
A blow-up function $\varphi$ is said to be polynomially bounded if the following inequality is satisfied for all $t\in[0,T)$ with $p_1$ and $p_2$ being polynomials:
\begin{align}
0 < p_1\left(\frac{T}{T-t}\right) \leq \varphi(t) \leq p_2\left(\frac{T}{T-t}\right).
\end{align}
The time scale transformation $\tau = a(t)$ is said to be a polynomially bounded time scale transformation if its derivative $a'(t)$  is a polynomially bounded blow-up function.
\end{definition}

\begin{definition}
  \label{defn:polynomially_bounded_semi}
  Polynomially Bounded Semi-Blow-Up Function:
  A function $\varphi$ is said to be a polynomially bounded semi-blow-up
  function if the inequality
  $|\varphi(t)| \leq \overline\varphi(t)$
  is satisfied for all $t\in[0,T)$ with
  $\overline\varphi$ being a polynomially bounded blow-up function.
  This implies that a polynomial $p_2$ exists such that $|\varphi(t)|\leq
  p_2\left(\frac{T}{T-t}\right)$ for all $t\in[0,T)$.
\end{definition}
Note that all polynomially bounded blow-up functions are also polynomially bounded semi-blow-up functions (but not necessarily vice versa).

While prescribed-time convergence was introduced in Definition~\ref{defn:ptc}, a more stringent ({\em exponential}) requirement for convergence introduced in the definition below will be seen to be crucial when considering interconnections of prescribed-time convergent systems. 

\begin{definition}
  \label{defn:prescribed_time_exp_convergence}
  Prescribed-Time Exponentially Convergent Signal:
A signal $q(t)$ is said to be prescribed-time exponentially convergent (to zero) if the inequality
  $|q(t)| \leq c e^{-\int_{\tilde T}^t \varphi(s)ds}$
holds for all $t\in(\tilde T,T)$ with some $0 \leq \tilde T < T$ and $c$ being a
positive constant and with $\varphi$ being a blow-up function such that over the
time interval $(\tilde T,T)$, the inequality 
$\varphi(t) \geq p_0\left(\frac{T}{T-t}\right)^2$
is satisfied with a positive constant $p_0$.
\end{definition}

It is trivially seen that if $x(t)$ in system \eqref{eq:sys_without_exogenous}
is a prescribed-time exponentially convergent signal, then the system
\eqref{eq:sys_without_exogenous} is PT-C. However, the reverse does not follow
in general. It is interesting to note that if $q(t)$ is a prescribed-time
exponentially convergent signal, then from the viewpoint of the time scale
transformation constructed as $\tau=a(t)=\int_0^t\varphi(s)ds$, we can write the inequality
  $|\tilde q(\tau)| \leq c e^{a(\tilde T)}e^{-\tau}$
where $\tilde q(\tau)=q(a^{-1}(\tau))$ is the signal $q$ expressed in the new
time variable $\tau$. This inequality in terms of $\tilde q(\tau)$ is simply the standard exponential convergence property over the time interval $[0,\infty)$ in terms of the new time variable $\tau$.

\begin{lemma}
  \label{lemma:polynomial_exp_convergence}
  The product of a polynomially bounded semi-blow-up function $\varphi(t)$ and a prescribed-time exponentially convergent signal $q(t)$ goes to 0 as $t\rightarrow T$, i.e., $\lim_{t\rightarrow T}\varphi(t)q(t)=0$. 
  
  \noindent{\bf Proof of Lemma~\ref{lemma:polynomial_exp_convergence}:} By the definition of a prescribed-time exponentially convergent signal, we know that the inequality $|q(t)|\leq c e^{-\int_{\tilde T}^t \tilde\varphi(s)ds}$ is satisfied with constants $c>0$ and $\tilde T\in[0,T)$ and with $\tilde\varphi(t)$ being a blow-up function such that $\tilde\varphi(t)\geq \tilde p_0\xi^2$ where $\xi\eqdef \frac{T}{T-t}$ and $\tilde p_0$ is a positive constant. Note that the interval $t\in[0,T)$ maps to $\xi\in[1,\infty)$. Since $\varphi(t)$ is a polynomially bounded semi-blow-up function, we know that $|\varphi(t)|\leq p_2(\xi)$ with $p_2$ being a polynomial. Defining the change of variables $\rho=\frac{T}{T-s}$, note that we have $\tilde\varphi(s)\geq \tilde p_0\rho^2$. Since $d\rho=\frac{\rho^2}{T}ds$, we obtain $\int_{\tilde T}^t \tilde\varphi(s)ds\leq\int_{\tilde T_1}^\xi \tilde p_0 T d\rho=\tilde p_0 T(\xi-\tilde T_1)$ where $\tilde T_1=\frac{T}{T-\tilde T}$. Note that $\tilde T_1\in[1,\infty)$. Now, we have $|\varphi(t)q(t)|\leq p_2(\xi)e^{-\tilde p_0 T(\xi-\tilde T_1)}$. Hence, $|\varphi(t)q(t)|\rightarrow 0$ as $\xi\rightarrow \infty$, i.e., as $t\rightarrow T$. 
\end{lemma}

Indeed, a stronger statement
% than Lemma~\ref{lemma:polynomial_exp_convergence}
holds as shown below.

\begin{lemma}
  \label{lemma:polynomial_exp_convergence2}
  The product of a polynomially bounded semi-blow-up function $\varphi(t)$ and a prescribed-time exponentially convergent signal $q(t)$ is a prescribed-time exponentially convergent signal.

  \noindent{\bf Proof of Lemma~\ref{lemma:polynomial_exp_convergence2}:} Using
  the notation in the proof of Lemma~\ref{lemma:polynomial_exp_convergence}, we have $|q(t)|\leq c e^{-\int_{\tilde T}^t \tilde\varphi(s)ds}$ over the time interval $t\in(\tilde T,T)$. Define $\tilde q(t)=\frac{q(t)}{e^{-0.5\int_{\tilde T}^t \tilde\varphi(s)ds}}$. Then, we have $|\tilde q(t)|\leq c e^{-0.5\int_{\tilde T}^t \tilde\varphi(s)ds}$ over the time interval $t\in(\tilde T,T)$. Hence, $\tilde q(t)$ is a prescribed-time exponentially convergent signal implying from Lemma~\ref{lemma:polynomial_exp_convergence} that $\lim_{t\rightarrow T}\varphi(t)\tilde q(t)=0$. Therefore, after some time $\tilde T_1$, we have $|\varphi(t)\tilde q(t)|\leq 1$ for all $t\in[\tilde T_1,T)$. Hence, we have $|\varphi(t)q(t)|=|\varphi(t)\tilde q(t)|e^{-0.5\int_{\tilde T}^t \tilde\varphi(s)ds} \leq e^{-0.5\int_{\tilde T}^t \tilde\varphi(s)ds}$ for all $t\in[\max\{\tilde T,\tilde T_1\},T)$. Therefore, $\varphi(t)q(t)$ is a prescribed-time exponentially convergent signal.
\end{lemma}

\begin{lemma}
  \label{lemma:integral_polynomial_boundedness}
  If $\varphi(t)$ is a polynomially bounded semi-blow-up function and $\nu(t):[0,T)\rightarrow {\cal R}$ is a uniformly bounded function, then $\int_0^t \varphi(s)\nu(s)ds$ is a polynomially bounded semi-blow-up function.

  \noindent{\bf Proof of Lemma~\ref{lemma:integral_polynomial_boundedness}:} Since $\varphi(t)$ is a polynomially bounded semi-blow-up function, we have $|\varphi(t)|\leq p_2(\xi)$ with $p_2$ being a polynomial and $\xi=\frac{T}{T-t}$. Since $\nu(t)$ is uniformly bounded, we have $|\int_0^t \varphi(s)\nu(s)ds|\leq \sup_{s\in[0,T)}|\nu(s)|\int_0^t |\varphi(s)|ds$. Defining the change of variables $\rho=\frac{T}{T-s}$, we have $d\rho=\frac{\rho^2}{T}ds$. Hence, we obtain $|\int_0^t \varphi(s)\nu(s)ds|\leq T\sup_{s\in[0,T)}|\nu(s)|\int_1^\xi \frac{p_2(\rho)}{\rho^2} d\rho$.
  Note that if the polynomial $p_2(\rho)$ is of the form $\sum_{i=0}^Na_i
  \rho^i$, then $\int_1^\xi \frac{p_2(\rho)}{\rho^2} d\rho\leq
  a_0+(a_1+a_2)\xi+\sum_{i=3}^N \frac{a_i}{i-1}\xi^{i-1}$. Hence, we see that $\int_0^t \varphi(s)\nu(s)ds$ is a polynomially bounded semi-blow-up function.
\end{lemma}

\noindent{\bf Remark 1:} A class of blow-up functions that has been widely used in the literature is of the form $\varphi(t)=c\Big(\frac{T}{T-t}\Big)^k$ with constant $k \geq 1$ and constant $c>0$.
From Definition~\ref{defn:polynomially_bounded}, it is seen that these blow-up functions are polynomially bounded.
The corresponding time scale transformations defined as $a(t)=\int_0^t \varphi(s)ds$
are of form $a(t)=-cT\log\Big(\frac{T-t}{T}\Big)$ when $k=1$ and
$a(t)=\frac{c}{k-1}\Big[T\Big(\frac{T}{T-t}\Big)^{k-1}-T\Big]$ when $k>1$.
Defining $\tau=a(t)$, the inverse functions are seen to be
$t=T(1-e^{-\frac{\tau}{cT}})$ for $k=1$
and $t=T\Big\{1-\Big(\frac{cT}{(k-1)\tau+cT}\Big)^{\frac{1}{k-1}}\Big\}$ for $k>1$.
Hence, in terms of the transformed time variable $\tau$, we have $\varphi(a^{-1}(\tau))=e^{\frac{\tau}{cT}}$ for $k=1$
and $\varphi(a^{-1}(\tau))=\Big(\frac{(k-1)\tau+cT}{cT}\Big)^{\frac{k}{k-1}}$ for $k>1$.
It can be noted that the functions $\varphi(\tau)$ are also polynomially bounded
in $\tau$ when $k>1$. The property of a class of blow-up functions being
polynomially bounded in the new time variable $\tau$ is seen to be crucial in
control designs such as \cite{KKK19a,KKK20a} where polynomial boundedness along with
prescribed-time exponential convergence of a Lyapunov function were instrumental
in inferring prescribed-time convergence of the system state to zero.
Also, note that since $\frac{T}{T-t}\geq 1$ over $t\in[0,T)$, a blow-up function
of form $\varphi(t)=c\Big(\frac{T}{T-t}\Big)^k$ with constant $c > 0$ and
constant $k\geq 2$
also satisfies the inequality on $\varphi(t)$ introduced in the
definition of a prescribed-time exponentially convergent signal (Definition~\ref{defn:prescribed_time_exp_convergence}).

We now proceed to formulating Lyapunov-based characterizations for the prescribed-time properties summarized in Section~\ref{sec:preliminaries}. Unlike the corresponding analysis for asymptotic stability, Lyapunov functions in the context of prescribed-time analysis are often explicitly time-dependent and are not positive-definite functions of $x$ as discussed further in Remark 2. Hence, in the lemmas below, we first state more general properties in terms of non-negative functions $V$ and then specialize to the case where $V$ is a positive-definite function. 

\begin{lemma}
  \label{lemma:V_converge}
  Consider system \eqref{eq:sys_without_exogenous}. If a non-negative function $V(x,t)$ satisfies the following
  inequality along the system trajectories with $\varphi$ being a blow-up function satisfying $\varphi(t) \geq p_0\left(\frac{T}{T-t}\right)^2$ with $p_0$ being a positive constant, then $V(x(t),t)$ is prescribed-time exponentially convergent to zero:
  \begin{align}
    \dot V \leq -\varphi(t)V.
    \label{eq:ptc_V}
  \end{align}
  \noindent{\bf Proof of Lemma~\ref{lemma:V_converge}:} The lemma can be proved
  in two ways which provide two different conceptual viewpoints. Firstly, using
  the Gr\"onwall-Bellman inequality, it follows from \eqref{eq:ptc_V} that $V\leq V_0 e^{-\int_0^t\varphi(s)ds}$ where $V_0=V(x(0),0)$. Hence, from the definition of a blow-up function, $V\rightarrow 0$ as $t \rightarrow T$ with furthermore prescribed-time exponential convergence as seen from Definition~\ref{defn:prescribed_time_exp_convergence}. Secondly, consider the time scale transformation $\tau = a(t)= \int_0^t \varphi(s)ds$. Denoting $\tilde a'(\tau)=a'(a^{-1}(\tau))$, we have $d\tau=\tilde a'(\tau)dt$. Hence, \eqref{eq:ptc_V} yields $\frac{dV}{d\tau}\leq -\frac{\varphi(t)}{\tilde a'(\tau)}V$. Noting that $\tilde a'(\tau)=\varphi(t)$ by the definition $\tau=a(t)$, we have $\frac{dV}{d\tau}\leq -V$ implying exponential convergence of $V$ to 0 in terms of the transformed time variable $\tau$. Since $\tau\rightarrow T$ corresponds to $t\rightarrow T$, it is seen that $V$ converges to 0 as $t\rightarrow T$ and furthermore $V\leq V_0 e^{-\tau}$ from which prescribed-time exponential convergence can be inferred.
\end{lemma}

\begin{lemma}
  \label{lemma:V_converge_iss}
  Consider system \eqref{eq:sys_with_exogenous}. If a non-negative function $V(x,t)$ satisfies the following
  inequality along the system trajectories with $\varphi$ being a blow-up function and with $a > 0$ and $b$ being constants
  \begin{align}
    \dot V \leq \varphi(t)[-a V + b |d|],
    \label{eq:pt_iss_V}
  \end{align}
  then the following inequality is satisfied for all $t\in[0,T)$:
  \begin{align}
  V(x(t),t) \leq \beta(|V_0|,a(t))+\gamma(\sup_{s\in[0,t)}|d(s)|)
  \label{eq:V_ineq_lemma}
  \end{align}
  where $a(t)=\int_0^t \varphi(s)ds$, $V_0=V(x(0),0)$,
  $\beta$ is a class ${\cal KL}$ function, and $\gamma$ is a class ${\cal K}$ function.
  Furthermore, if $\lim_{t\rightarrow T}|d(t)|=0$, then $\lim_{t\rightarrow T}V(x(t),t)=0$.

  \noindent{\bf Proof of Lemma~\ref{lemma:V_converge_iss}:} As in the proof of Lemma~\ref{lemma:V_converge}, consider the time scale transformation $\tau = a(t)= \int_0^t \varphi(s)ds$ with which we obtain (analogously to the proof of Lemma~\ref{lemma:V_converge})
\begin{align}
  \frac{dV}{d\tau} \leq -a V +b |d|.
  \label{eq:Vdot_tau_PT_ISS}
\end{align}
Hence, in terms of the transformed time variable $\tau$, the system satisfies the standard ISS property over the infinite time interval $\tau\in[0,\infty)$. Therefore, we know that for all time instants $\tau$, an inequality of the form \eqref{eq:V_ineq_lemma} is satisfied with $\beta$ being a class ${\cal KL}$ function and $\gamma$ a class ${\cal K}$ function.
In addition, when $\lim_{t\rightarrow T}|d(t)|=0$, we see from the dynamics of $V$ in terms of the transformed time variable $\tau$ as given by \eqref{eq:Vdot_tau_PT_ISS} which is an asymptotically stable system with an exogenous input converging to 0 as $\tau\rightarrow\infty$ that therefore also $V$ converges to 0 as $\tau\rightarrow\infty$, i.e., as $t\rightarrow T$.
\end{lemma}

\begin{lemma}
  \label{lemma:V_converge_iss_c}
  Consider system \eqref{eq:sys_with_exogenous}. If a non-negative function $V(x,t)$ satisfies the following
  inequality along the system trajectories with $\varphi$ being a blow-up function and with $a > 0$ and $b$ being constants, then $\lim_{t\rightarrow T}V(x(t),t)=0$:
  \begin{align}
    \dot V \leq \varphi(t)[-a V] + b |d|.
    \label{eq:pt_iss_c_V}
  \end{align}
\noindent{\bf Proof of Lemma~\ref{lemma:V_converge_iss_c}:} Defining the time scale transformation $\tau = a(t)= \int_0^t \varphi(s)ds$ analogously to the proof of Lemma~\ref{lemma:V_converge}, we obtain
$\frac{dV}{d\tau} \leq -a V +\frac{b d}{\varphi(a^{-1}(\tau))}$. Noting that in
terms of the transformed time variable $\tau$, this corresponds to a stable linear system driven by an exogenous input that converges to 0 as $\tau\rightarrow\infty$, it follows that $V$ goes to 0 as $\tau\rightarrow \infty$, i.e., as $t\rightarrow T$.
\end{lemma}

\begin{lemma}
  \label{lemma:ptc}
  Consider system \eqref{eq:sys_without_exogenous}. Under the conditions of Lemma~\ref{lemma:V_converge}, if furthermore $V(x,t)$ is a positive-definite function\footnote{A function $V(x,t):{\cal R}^n\times [0,T)\rightarrow [0,\infty)$ with $t$ being time is said to be a positive-definite function if $V(0,t)=0$ and $V(x,t)\geq \alpha(|x|)$ for all $x\in{\cal R}^n$ and $t\in [0,T)$ with
  $\alpha$ being a class ${\cal K}_\infty$ function.}, then the system
  \eqref{eq:sys_without_exogenous} is PT-C.

  \noindent{\bf Proof of Lemma~\ref{lemma:ptc}:} From Lemma~\ref{lemma:V_converge}, we see that $V\rightarrow 0$ as $t\rightarrow T$. Since $V$ is a positive-definite function, we see that the PT-C property follows.
\end{lemma}

\begin{lemma}
  \label{lemma:pt_iss}
  Consider system \eqref{eq:sys_with_exogenous}. Under the conditions of Lemma~\ref{lemma:V_converge_iss}, if furthermore $V(x,t)$ is a positive-definite function, then the system
  \eqref{eq:sys_with_exogenous} is PT-ISS.
  Also, if $\lim_{t\rightarrow T}|d(t)|=0$, then $\lim_{t\rightarrow T}|x(t)|=0$.

  \noindent{\bf Proof of Lemma~\ref{lemma:pt_iss}:}
This lemma follows from Lemma~\ref{lemma:V_converge_iss} by noting that if $V(x,t)$
is a positive-definite function, then the inequality \eqref{eq:Vdot_tau_PT_ISS} implies that the PT-ISS property in Definition~\ref{defn:pt_iss} is satisfied with blow-up function $\mu(t)=a(t)+1$.
\end{lemma}

\begin{lemma}
  \label{lemma:pt_iss_c}
  Consider system \eqref{eq:sys_with_exogenous}. Under the conditions of Lemma~\ref{lemma:V_converge_iss_c}, if furthermore $V(x,t)$ is a positive-definite function, then the system
  \eqref{eq:sys_with_exogenous} is PT-ISS+C.
  
\noindent{\bf Proof of Lemma~\ref{lemma:pt_iss_c}:} Using Lemma~\ref{lemma:V_converge_iss_c}, it is seen that if $V(x,t)$ is a positive-definite function, the PT-ISS+C property is satisfied.
\end{lemma}

\noindent{\bf Remark 2:} Analogous to the use of Lyapunov functions in demonstrating asymptotic stability properties (as $t\rightarrow\infty$) of nonlinear systems, it can be expected that Lyapunov (or Lyapunov-like) functions would play a useful role to demonstrate the various prescribed-time stability and convergence properties summarized in the definitions above. However, while typical Lyapunov functions that appear in asymptotic stability analysis often do not involve the time $t$ explicitly, the Lyapunov functions that appear most naturally for prescribed-time analysis are often (but not always) time-varying.
To illustrate cases where Lyapunov functions for prescribed-time analysis do not and do include the time $t$ explicitly, we consider two examples. As a first example, consider the scalar system $\dot x = u$ and consider a control law of form $u=-\varphi(t)x$ with $\varphi(t)$ being a blow-up function as in Lemma~\ref{lemma:V_converge}. Defining the Lyapunov function $V=\frac{1}{2}x^2$, we obtain $\dot V = -\varphi(t)x^2 = -2\varphi(t)V$ which is of the form considered in Lemma~\ref{lemma:V_converge} implying that $V\rightarrow 0$ as $t\rightarrow T$ and therefore also 
$x\rightarrow 0$ as $t\rightarrow T$. 
As a second example, consider the second-order system $\dot x_1=x_2, \dot x_2=u$. A control law can be designed using the backstepping technique \cite{KKK95,JK95,KK03} by first defining the virtual control law $x_2^*=-\varphi(t)x_1$ with $\varphi(t)$ being a blow-up function. Then, define $z_2=x_2-x_2^*$ and define the Lyapunov function 
\begin{align}
V &=
\frac{1}{2}x_1^2+\frac{c}{2}z_2^2
= 
\frac{1}{2}x_1^2+\frac{c}{2}(x_2+\varphi(t)x_1)^2
\end{align}
with $c$ being a positive constant. Then, we obtain 
\begin{align}
\dot V &= -\varphi(t)x_1^2 + x_1z_2 + cz_2[u+\varphi(t)x_2+\dot\varphi(t)x_1].
\end{align}
Designing the control law $u=-\frac{x_1}{c}-\varphi(t)x_2 - \dot\varphi(t)x_1 - \varphi(t)z_2$, we obtain 
$\dot V = -2\varphi(t)V$ which implies from Lemma~\ref{lemma:V_converge} that $V\rightarrow 0$ as $t\rightarrow T$. However, unlike in the first example, it is not obvious that this implies that $x=[x_1,x_2]^T$ goes to 0 as $t\rightarrow T$ since $V$ is not only a function of $x$, but also explicitly involves the time $t$ via the term $\varphi(t)$ that appears in the definition of $V$. Furthermore, although $V$ evidently can be written as a quadratic form in $x$ as
\begin{align}
  V &= x^T \overline P(t) x \ \ , \ \ \overline P(t) = \frac{1}{2} \begin{bmatrix}
    1+c\varphi^2(t) & c\varphi(t) \\ c\varphi(t) & c
  \end{bmatrix},
\end{align}
it can be seen that no positive-definite matrix $P$ exists such that $V\geq x^T Px$. If such a positive-definite matrix $P$ existed, then the convergence of $V$ to 0 would immediately imply the convergence of $x$ to 0. To see that such a positive-definite matrix $P$ can not exist, note that the eigenvalues of the matrix $\overline P(t)$ are given by the roots of its characteristic polynomial $\lambda^2 - [c+1+c\varphi^2(t)]\lambda + c$. It can be easily shown that one root of this characteristic polynomial (i.e., one eigenvalue of $\overline P(t)$) goes to 0 as $\varphi(t)\rightarrow \infty$. Hence, a constant positive-definite matrix $P$ does not exist such that $V\geq x^TPx$. However, it can still be inferred that $x\rightarrow 0$ as $t\rightarrow T$ by using the fact that $V\rightarrow 0$ as $t\rightarrow T$ \underline{\em and} the fact that the convergence of $V$ to 0 is {\em fast} in the sense of the following inequality that follows from $\dot V = -2\varphi(t)V$:
 \begin{align}
   V_t\leq V_0 e^{-2\int_0^t\varphi(s)ds}
   \label{eq:remark_V_converge}
 \end{align}
 where $V_0$ and $V_t$ denote the values of $V$ at time instants $0$ and $t$, respectively. To see this, note that from the convergence of $V$ to 0, we can directly infer the convergence of both $x_1$ and $x_2+\varphi(t)x_1$ to 0. Furthermore, since $\frac{1}{2}x_1^2 \leq V$ by the definition of $V$, we have $\frac{1}{2}x_1^2 \leq V_0 e^{-2\int_0^t\varphi(s)ds}$, implying that $x_1$ is prescribed-time exponentially convergent to 0. If $\varphi(t)$ is picked to be a polynomially bounded blow-up function, we then see from Lemma~\ref{lemma:polynomial_exp_convergence} that $\varphi(t)x_1\rightarrow 0$ as $t\rightarrow T$. Therefore, it follows also that $x_2\rightarrow 0$ as $t\rightarrow T$.  The inequality \eqref{eq:remark_V_converge} implies that $V$ is prescribed-time exponentially convergent to 0 by Definition~\ref{defn:prescribed_time_exp_convergence}. From the above analysis, it is seen that for this example, this {\em fast} convergence of $V$ is crucial to infer prescribed-time convergence of $x$ to 0 while just the fact that $V$ converges to 0 is not sufficient to infer the convergence property of $x$. This notion of a Lyapunov function providing a convergence certificate when the convergence of $V$ is fast enough is formalized in the definition of a prescribed-time exponentially convergent Lyapunov certificate below.

 \begin{definition}
   \label{defn:prescribed_time_exp_convergent_certificate}
   Prescribed-Time Exponentially Convergent Lyapunov Certificate:
   A function $V(x(t),t)$ is said to be a prescribed-time exponentially convergent Lyapunov certificate for a system (e.g., of form \eqref{eq:sys_without_exogenous} or \eqref{eq:sys_with_exogenous}) if the fact that $V$ is a signal that is prescribed-time exponentially convergent (to zero) implies based on the form of $V$ and the system dynamics that $x\rightarrow 0$ as $t\rightarrow T$.
 \end{definition}

 \noindent{\bf Remark 3:} From Remark 2 and Definition~\ref{defn:prescribed_time_exp_convergent_certificate}, it is seen that the $V$ constructed in the second example in Remark 2 is an example of a prescribed-time exponentially convergent Lyapunov function. Note that a positive-definite function $V$ is (trivially) a prescribed-time exponentially convergent Lyapunov function.
 It is also to be noted that while not explicitly stated in the definition of a prescribed-time exponentially convergent Lyapunov function, the crucial element of such a function in the second example in Remark 2 is that the blow-up function which appears in $V$ and introduces the time dependence of $V$ is polynomially bounded. This polynomial boundedness enables inferring of asymptotic convergence of the system state from the prescribed-time exponential convergence of $V$ (essentially the fact that the product of a polynomially bounded increasing function and an exponentially decreasing function tends to zero, i.e., Lemma~\ref{lemma:polynomial_exp_convergence}). This is, in spirit, analogous to the prescribed-time high-gain based control designs in \cite{KKK19a,KKK20a} where the dynamic high-gain scaling parameter (which can be viewed as playing an analogous role to a blow-up function) was shown to be polynomially bounded in the new time variable $\tau$ and it was shown that a Lyapunov function (involving the high-gain scaling parameter) converges exponentially in terms of $\tau$. This pair of properties (polynomial boundedness of dynamic high-gain scaling parameter, exponential convergence of Lyapunov function) was crucial in \cite{KKK19a,KKK20a} to infer prescribed-time convergence of the system state to zero. While both these properties were formulated in terms of the new time variable $\tau$ in \cite{KKK19a,KKK20a}, the analogous properties introduced in Definitions~\ref{defn:polynomially_bounded} and \ref{defn:prescribed_time_exp_convergence} do not involve a new time variable and are defined purely in terms of $t$. This is crucial for considering interconnected systems in this paper since the different subsystems in the interconnected combination involve, in general, different blow-up functions and therefore different ``intrinsic'' time scale transformations implying that there is no single new time variable $\tau$ that could be utilized.

Note that from Lemma~\ref{lemma:V_converge} and Definition~\ref{defn:prescribed_time_exp_convergent_certificate}, the following variant of Lemma~\ref{lemma:ptc} follows.
 
\begin{lemma}
   \label{lemma:ptc_exp}
   Consider system \eqref{eq:sys_without_exogenous}. Under the conditions of Lemma~\ref{lemma:V_converge}, if furthermore $V(x,t)$ is a prescribed-time exponentially convergent Lyapunov certificate for the system \eqref{eq:sys_without_exogenous}, then the system
   \eqref{eq:sys_without_exogenous} is PT-C.

   \noindent{\bf Proof of Lemma~\ref{lemma:ptc_exp}:} From Lemma~\ref{lemma:V_converge}, we see that $V$ is prescribed-time exponentially convergent to zero, therefore implying the PT-C property of \eqref{eq:sys_without_exogenous} since $V$ is given to be a prescribed-time exponentially convergent Lyapunov certificate for the system \eqref{eq:sys_without_exogenous}.
\end{lemma}

Also, the following variant of Lemma~\ref{lemma:pt_iss} holds. Note that while it was sufficient in Lemma~\ref{lemma:pt_iss} that $d$ converges to zero as $t\rightarrow T$ to be able to infer that the state $x$ goes to zero, the lemma below requires the stronger condition that $d$ is prescribed-time {\em exponentially} convergent to zero as $t\rightarrow T$.

\begin{lemma}
  \label{lemma:ptc_exp_with_d}
  Consider system \eqref{eq:sys_with_exogenous}. Under the conditions of Lemma~\ref{lemma:V_converge_iss}, if furthermore $\varphi$ is a polynomially bounded blow-up function satisfying $\varphi(t) \geq p_0\left(\frac{T}{T-t}\right)^2$ with $p_0$ being a positive constant, $V(x,t)$ is a prescribed-time exponentially convergent Lyapunov certificate for the system \eqref{eq:sys_with_exogenous}, and $d$ is  uniformly bounded over $[0,T)$ and a prescribed-time exponentially convergent signal, then along all trajectories of the system
  \eqref{eq:sys_with_exogenous}, we have $\lim_{t\rightarrow T}|x(t)|=0$.

  \noindent{\bf Proof of Lemma~\ref{lemma:ptc_exp_with_d}}: By the
  Gr\"onwall-Bellman inequality\cite{Bellman53,Pachpatte98}, \eqref{eq:pt_iss_V} implies that 
  $V\leq [V_0+\int_0^t\varphi(s)|b d(s)|ds]e^{-\int_0^t a \varphi(s)ds}$ where $V_0=V(x(0),0)$. By Lemma~\ref{lemma:integral_polynomial_boundedness}, we see that $\nu(t)\eqdef\int_0^t\varphi(s)|b d(s)|ds$ is a polynomially bounded semi-blow-up function since $\varphi(s)$ is a polynomially bounded blow-up function and $|bd(s)|$ is a uniformly bounded signal. By Lemma~\ref{lemma:polynomial_exp_convergence2}, we then see that $\nu(t)e^{-\int_0^t a \varphi(s)ds}$ is therefore a prescribed-time exponentially convergent signal. Since $V(x,t)$ is given to be a prescribed-time exponentially convergent Lyapunov certificate for the system \eqref{eq:sys_with_exogenous}, this implies that $x\rightarrow 0$ as $t\rightarrow T$.
\end{lemma} 

\noindent{\bf Remark 4:} While the prescribed-time convergence and ISS properties above have been worded with the entire state $x$ as the quantity regarding which the prescribed-time properties are to be attained, it is to be noted that, in general, it could be a subset of the state vector that is relevant for enforcement of these properties. An example of such a scenario is wherein the prescribed-time properties are achieved using a control law that utilizes a dynamic gain which would be a controller state variable. For example, in our dynamic scaling-based control designs in \cite{KKK19a,KKK20a}, the control law utilizes powers of a dynamic scaling parameter $r$ whose dynamics is constructed as part of the design of the dynamic control law. By design of the scaling parameter dynamics, $r$ goes to $\infty$ as $t\rightarrow T$. However, the dynamics of $r$ are designed such that $r$ grows at most polynomially in terms of a transformed time variable $\tau=a(t)$ while the convergence of scaled state variables (defined as products of state variables and powers of $r$) is shown to be exponential in terms of $\tau$. Hence, by appropriately controlling the growth rate of $r$, it is shown in the control design approaches in \cite{KKK19a,KKK20a} that the original system state and the control law go to zero as $t\rightarrow T$ even though $r\rightarrow \infty$ as $t\rightarrow T$. In scenarios such as this, the prescribed-time convergence and ISS properties would only consider the original system's state and not the state variable $r$ which is part of the dynamic control law.  % For brevity and clarity, we will use the terminology ``core state'' to refer to the part of the overall state vector that is to be considered in the prescribed-time properties.   

\noindent{\bf Remark 5:} For simplicity and brevity, the forms of the Lyapunov inequalities in Lemmas~\ref{lemma:V_converge}--\ref{lemma:ptc_exp_with_d} were written with $V$ appearing linearly on the right hand sides of the Lyapunov inequalities. However, these Lyapunov inequalities can be generalized to instead have nonlinear functions of $V$ on the right hand sides albeit with more algebraic complexity in the proofs of the lemmas.

\vspace*{-0.1in}
\section{Cascade and Feedback Interconnections of PT-ISS Systems}
\label{sec:interconnections}
In this section, we develop the main results of this paper. We will first consider two interconnected systems of the following forms and will then consider the more general case of an arbitrary number of interconnected systems.
\begin{itemize}
\item cascade interconnection of a PT-C system and a PT-ISS system:
  \begin{align}
    \dot x_1 &= f_1(x_1,t) \ \ \ ; \ \ \ 
    \dot x_2 = f_2(x_2,x_1,t).
               \label{eq:cascade_interconnection_12}
  \end{align}
\item feedback interconnection of two PT-ISS systems:
  \begin{align}
    \dot x_1 &= f_1(x_1,x_2,t) \ \ \ ; \ \ \ 
    \dot x_2 = f_2(x_2,x_1,t).
               \label{eq:feedback_interconnection_12}
  \end{align}
\end{itemize}
In both the cases above, $x_1\in{\cal R}^{n_1}$ and $x_2\in{\cal R}^{n_2}$ are the state vectors of the two systems.
For each of the two forms of interconnected systems above, we will first consider the simpler case where the blow-up functions corresponding to the two systems are the same and then the general case where the blow-up functions are different.

\begin{thrm}
  \label{thm:cascade_same_blowup}
  Consider the cascade interconnection shown in
  \eqref{eq:cascade_interconnection_12}.
  If prescribed-time exponentially convergent Lyapunov certificates $V_1(x_1,t)$
  and $V_2(x_2,t)$, respectively, for the two systems exist satisfying
  \begin{align}
    \dot V_1 &\leq -\varphi(t) V_1 \ \ \ ; \ \ \ 
    \dot V_2 \leq \varphi(t)[-a V_2 + b V_1]
               \label{eq:cascade_V12}
  \end{align}
  with $\varphi$ being a blow-up function
  satisfying $\varphi(t) \geq p_0\left(\frac{T}{T-t}\right)^2$ with $p_0$ being a positive constant, and with
  $a>0$ and $b$ being constants, then the overall system formed by the cascade interconnection is PT-C.

\noindent{\bf Proof of Theorem~\ref{thm:cascade_same_blowup}:} Defining
$V=c_1V_1+c_2V_2$ with $c_1$ and $c_2$ being positive constants such that
$\frac{c_1}{c_2}>b$, we see that $\dot V \leq -\kappa \varphi(t)V$ where
$\kappa=\min\Big\{\frac{c_1-c_2b}{c_1},a\Big\}$ is a positive constant by the
construction of $c_1$ and $c_2$, from which it can be inferred using Lemma~\ref{lemma:V_converge} that
$V$ is a prescribed-time exponentially convergent signal. From the definition of $V$, $V_1$ and $V_2$ are also prescribed-time exponentially convergent signals. Since $V_1$ and $V_2$ are prescribed-time exponentially convergent Lyapunov certificates for the first and second system, respectively, it follows that each of the systems is PT-C and therefore so is the overall interconnected system.
\end{thrm}

In the theorem below, we note that a stronger statement than
Theorem~\ref{thm:cascade_same_blowup} holds in which the blow-up functions
appearing in the time derivatives of the two Lyapunov functions are allowed to
be different and a nonlinear function of $V_1$ is allowed to appear in $\dot
V_2$. However, the blow-up function appearing in $\dot V_2$ is required to be
polynomially bounded. 

\begin{thrm}
  \label{thm:cascade_different_blowup}
  Consider the cascade interconnection shown in \eqref{eq:cascade_interconnection_12}. If prescribed-time exponentially convergent Lyapunov certificates $V_1(x_1,t)$ and $V_2(x_2,t)$, respectively, for the two systems exist satisfying
  \begin{align}
    \dot V_1 &\leq -\varphi_1(t) V_1 \ \ \ ; \ \ \ 
    \dot V_2 \leq \varphi_2(t)[-a V_2 + \varphi_3(t)p(V_1)]
               \label{eq:cascade_different_V12}
  \end{align}
  with $a>0$ being a constant, $p(.)$ a polynomial such that $p(0)=0$, $\varphi_1$ a blow-up function and $\varphi_2$ a polynomially
  bounded blow-up function
  satisfying $\varphi_1(t) \geq p_{01}\left(\frac{T}{T-t}\right)^2$ and
  $\varphi_2(t) \geq p_{02}\left(\frac{T}{T-t}\right)^2$
  with $p_{01}$ and $p_{02}$ being positive constants,
  and with $\varphi_3$
  being a polynomially bounded semi-blow-up function, then the overall system formed by the cascade interconnection is PT-C.

\noindent{\bf Proof of Theorem~\ref{thm:cascade_different_blowup}:} From
\eqref{eq:cascade_different_V12}, we know that $V_1$ is uniformly bounded and we
also know using Lemma~\ref{lemma:V_converge} that $V_1$ is a prescribed-time
exponentially convergent signal. Since $p(.)$ is a polynomial with $p(0)=0$, it is evident
that $p(V_1)$ is also a prescribed-time
exponentially convergent signal. Using Lemma~\ref{lemma:polynomial_exp_convergence2}, we see that $\varphi_3(t)p(V_1)$ is also a prescribed-time
exponentially convergent signal since $\varphi_3(t)$ is a polynomially bounded
semi-blow-up function. Hence, using Lemma~\ref{lemma:ptc_exp_with_d}
and noting that $\varphi_2(t)$ is a polynomially bounded blow-up function, we
see that $V_2$ is also a prescribed-time exponentially convergent signal. Since
$V_1$ and $V_2$ are prescribed-time exponentially convergent Lyapunov certicates
for the two systems, it follows that the two systems (and therefore the overall
interconnected system) are PT-C\footnote{Note however that unlike the proof of Theorem~\ref{thm:cascade_same_blowup}, exhibiting a Lyapunov function for the overall interconnected system can not be accomplished with just a linear constant-coefficient combination, but would need a more general construction analogous to Theorem~\ref{thm:feedback_different_blowup}.}. 
\end{thrm}

% While Theorems~\ref{thm:cascade_same_blowup} and
% \ref{thm:cascade_different_blowup} considered cascade interconnections, the
% following theorems consider properties of feedback interconnections of PT-ISS systems.

\begin{thrm}
  \label{thm:feedback_same_blowup}
  Consider the feedback interconnection shown in \eqref{eq:feedback_interconnection_12}. If prescribed-time exponentially convergent Lyapunov certificates $V_1(x_1,t)$ and $V_2(x_2,t)$, respectively, for the two systems exist satisfying
  \begin{align}
    \dot V_1 &\leq \varphi(t) [-a_1V_1 + b_1 V_2]
    \label{eq:feedback_V1}\\
    \dot V_2 &\leq \varphi(t)[-a_2 V_2 + b_2 V_1]
               \label{eq:feedback_V2}
  \end{align}
  with $\varphi$ being a blow-up function
  satisfying $\varphi(t) \geq p_0\left(\frac{T}{T-t}\right)^2$ with $p_0$ being a positive constant,
  and with $a_1$, $a_2$, $b_1$, and $b_2$ being positive constants\footnote{Without loss of
    generality, we consider $b_1$ and $b_2$ to be positive constants in Theorem~\ref{thm:feedback_same_blowup} since if, for example,
    $b_1\leq 0$, then \eqref{eq:feedback_V1} implies $\dot V_1\leq -a_1\varphi(t)V_1$ thus reducing the scenario considered in this theorem to the scenario in Theorem~\ref{thm:cascade_same_blowup} and similarly for the case $b_2\leq 0$. Analogously, we consider $b_1$ and $b_2$ in
    Theorem~\ref{thm:feedback_different_blowup} and $b_{i,j}$ in
    Theorems~\ref{thm:N_feedback_same_blowup} and \ref{thm:N_feedback_different_blowup} to be positive constants without loss of generality.}, then the overall system formed by the feedback interconnection is PT-C if $a_1a_2 > b_1b_2$.

\noindent{\bf Proof of Theorem~\ref{thm:feedback_same_blowup}:} Define
$V=c_1V_1+c_2V_2$ with $c_1$ and $c_2$ being positive constants such that
$\frac{c_1}{c_2}>\frac{b_2}{a_1}$
and
$\frac{c_2}{c_1}>\frac{b_1}{a_2}$.
Such a choice of $c_1$ and $c_2$ is possible due to the condition $a_1a_2 >
b_1b_2$ in the theorem.
We have $\dot V \leq
-\kappa\varphi(t)V$ where
$\kappa=\min\Big\{\frac{a_1c_1-b_2c_2}{c_1},\frac{a_2c_2-b_1c_1}{c_2}\}$.
Note that $\kappa$ is a positive constant by the choice of constants $c_1$
and $c_2$.
Hence,
using Lemma~\ref{lemma:V_converge}, we see that
$V$ is a prescribed-time exponentially convergent signal implying from the definition of $V$, that $V_1$ and $V_2$ are also prescribed-time exponentially convergent signals. Since $V_1$ and $V_2$ are prescribed-time exponentially convergent Lyapunov certificates for the first and second system, respectively, each of the two systems is PT-C and therefore the overall interconnected system is also PT-C.
\end{thrm}

\begin{thrm}
  \label{thm:feedback_different_blowup}
  Consider the feedback interconnection shown in \eqref{eq:feedback_interconnection_12}. If prescribed-time exponentially convergent Lyapunov certificates $V_1(x_1,t)$ and $V_2(x_2,t)$, respectively, for the two systems exist satisfying
  \begin{align}
    \dot V_1 &\leq \varphi_1(t) [-a_1V_1 + b_1 V_2]
    \label{eq:feedback_different_V1}\\
    \dot V_2 &\leq \varphi_2(t)[-a_2 V_2 + b_2 V_1]
               \label{eq:feedback_different_V2}
  \end{align}
  with $\varphi_1$ and $\varphi_2$ being polynomially bounded blow-up functions
  satisfying $\varphi_1(t) \geq p_{01}\left(\frac{T}{T-t}\right)^2$ and $\varphi_2(t) \geq p_{02}\left(\frac{T}{T-t}\right)^2$
  with $p_{01}$ and $p_{02}$ being positive constants, 
  and with $a_1$, $a_2$, $b_1$, and $b_2$ being positive constants, then the overall system formed by the feedback interconnection is PT-C if the following condition holds:
    $a_1a_2 > b_1b_2$.
    % \\
    % \lim_{t\rightarrow T}\varphi_1(t)e^{-\delta(A)\int_0^t\underline\varphi(s)ds} = 0
    % \label{eq:combination_bnd1_feedback_different_blowup}
    % \\
    % \lim_{t\rightarrow T}\varphi_2(t)e^{-\delta(A)\int_0^t\underline\varphi(s)ds} = 0
    % \label{eq:combination_bnd2_feedback_different_blowup}

\noindent{\bf Proof of Theorem~\ref{thm:feedback_different_blowup}:}
Let $A=\begin{bmatrix}
  -a_1 & b_1 \\ b_2 & -a_2
\end{bmatrix}$.
It can be seen that $-A$ is a Lyapunov diagonally stable matrix since if\footnote{ With $\eta_1,\ldots,\eta_k$ being real numbers, $\mbox{diag}(\eta_1,\ldots,\eta_k)$ denotes the $k\times k$ diagonal matrix with the $i^{th}$ diagonal element being $\eta_i$. } $P=\mbox{diag}(c_1,c_2)$ with $c_1$ and $c_2$ being positive constants, then $-(PA+A^TP) > 0$ if $\gamma=\sqrt{\frac{c_1}{c_2}}$ satisfies the equation $\gamma^2b_1+b_2<2\gamma\sqrt{a_1a_2}$. This equation admits a solution for $\gamma$ as a positive constant since $a_1a_2 > b_1b_2$. 
Define
\begin{align}
  \overline V &= \left[
                \frac{V_1}{\varphi_1(t)} , \frac{V_2}{\varphi_2(t)}\right]^T \ , \ V_c = [V_1,V_2]^T.
                \label{eq:Vdefn_feedback}
\end{align}
Note that $V_c=\chi \overline V$ where $\chi=\begin{bmatrix}
  \varphi_1(t) & 0 \\ 0 & \varphi_2(t)
\end{bmatrix}
$.
Since $\frac{d}{dt}\varphi_1(t)$ and $\frac{d}{dt}\varphi_2(t)$ are
non-negative for all time because $\varphi_1(t)$ and $\varphi_2(t)$ are
monotonically increasing by the definition of a blow-up function, we see
that
\begin{align}
  \dot{\overline V} &\leq_e AV_c
  \ \ \Longrightarrow \ \ 
  \dot{\overline V} \leq_e A\chi \overline V.
\end{align}
Now, let $q=q_\delta(A)$. Then, defining $V=q\overline V$ and noting that $q >_e 0$, we have
$\dot V  = q\dot{\overline V}\leq qA\chi \overline V\leq -\delta(A)
q\chi\overline V\leq -\delta(A) q\underline\varphi(t)\overline V$
where $\underline\varphi(t)\eqdef \min\{\varphi_1(t),\varphi_2(t)\}$
since $\delta(A)>0$ and $\overline V\geq_e 0$.  Therefore, noting that $\delta(A)$ and $\underline\varphi(t)$ are scalars, we have $\dot V \leq -\delta(A)\underline\varphi(t)V$. Hence,
\begin{align}
  V\leq e^{-\delta(A)\int_0^t\underline\varphi(s)ds}V_0
  \label{eq:Vbnd_feedback_different_blowup}
\end{align}
where $V_0$ is the value of $V$ at time $t=0$.
Now, note that by the definition of $V$, we have $V\geq q_i\frac{V_i}{\varphi_i(t)}, i=1,2$, where $q_i$ is the $i^{th}$ element of $q$. Therefore,
\begin{align}
  V_i \leq \frac{1}{q_i} \varphi_i(t)e^{-\delta(A)\int_0^t\underline\varphi(s)ds}V_0 \ , \ i=1,2.
  \label{eq:Vibnd_feedback_different_blowup}
\end{align}
We know that $\underline\varphi(s)ds\geq \underline p_0\left(\frac{T}{T-t}\right)^2$ where $\underline p_0=\min\{p_{01},p_{02}\}$. Hence, the right hand side of the inequality \eqref{eq:Vibnd_feedback_different_blowup} is a product of a polynomially bounded blow-up function and a prescribed-time exponentially convergent signal. Therefore, by Lemma~\ref{lemma:polynomial_exp_convergence2}, $V_i$ is a prescribed-time exponentially convergent signal for $i=1,2$. 
Since $V_1$ and $V_2$ are
prescribed-time exponentially convergent Lyapunov certificates for the first and
second system, respectively, it follows that each of the two systems is PT-C and
therefore, the overall interconnected system is also PT-C.
\end{thrm}

% \noindent{\bf Proof of Theorem~\ref{thm:feedback_different_blowup_poly}:} This theorem follows from Theorem~\ref{thm:feedback_different_blowup} by noting that the conditions \eqref{eq:combination_bnd1_feedback_different_blowup} and \eqref{eq:combination_bnd2_feedback_different_blowup} are automatically satisfied if $\varphi_1(t)$ and $\varphi_2(t)$ are polynomially bounded blow-up functions. To see this, note that from the definition of a blow-up function, positive constants $\varphi_{01}$ and $\varphi_{02}$ exist such that $\varphi_1(t)\geq \varphi_{01}$ and $\varphi_2(t)\geq \varphi_{02}$ for all $t\rightarrow[0,T)$. Hence, with
% $\underline\varphi(t)\eqdef \min\{\varphi_1(t),\varphi_2(t)\}$, we have $\underline\varphi(t)\geq \underline\varphi_0$ where $\underline\varphi_0 = \min\{\varphi_{01},\varphi_{02}\}$ and
% $e^{-\delta(A)\int_0^t\underline\varphi(s)ds}\leq e^{-\delta(A)\underline\varphi_0 t}$. Therefore, $\varphi_1(t)e^{-\delta(A)\int_0^t\underline\varphi(s)ds}$ and
% $\varphi_2(t)e^{-\delta(A)\int_0^t\underline\varphi(s)ds}$ are products of polynomially bounded functions ($\varphi_1(t)$, $\varphi_2(t)$) and

While Theorems~\ref{thm:cascade_same_blowup}-\ref{thm:feedback_different_blowup} considered an interconnection of two PT-ISS systems, we now consider the general interconnection of an arbitrary number ($N$) of PT-ISS systems of form:
\begin{align}
  \dot x_i &= f_i (x_i,x_1,\ldots,x_{i-1},x_{i+1},\ldots,x_N,t)
             \label{eq:N_feedback_interconnection_i}
\end{align}
where the right hand side of the dynamics of state vector $x_i$ of the
$i^{th}$ system depends on $x_i$ as well as state vectors $x_j$, $j\neq i$,
of the other systems in the overall interconnected system.
% Then, the generalizations of Theorems~\ref{thm:feedback_same_blowup}--\ref{thm:feedback_different_blowup} are developed below.

\begin{thrm}
  \label{thm:N_feedback_same_blowup}
  Consider the feedback interconnection of $N$ systems with the dynamics of the
  $i^{th}$ system as shown in \eqref{eq:N_feedback_interconnection_i}. If
  prescribed-time exponentially convergent Lyapunov certificates $V_i(x_i,t)$
  for each of the systems ($i=1,\ldots,N$) exist satisfying
  \begin{align}
    \dot V_i &\leq \varphi(t) [-a_iV_i + \sum_{\substack{j=1\\j\neq i}}^{j=N} b_{i,j} V_j] \ , \ i=1,\ldots,N
    \label{eq:N_feedback_Vi}
  \end{align}
  with $\varphi$ being a blow-up function satisfying $\varphi(t) \geq
  p_0\left(\frac{T}{T-t}\right)^2$ with $p_0$ being a positive constant, and with
  $a_i, i=1,\ldots,N$ and $b_{i,j}, i=1,\ldots,N, j=1,\ldots,N, j\neq i$, being positive constants, then the overall system formed by the feedback interconnection is PT-C if the matrix $A$ defined below is a strict Hurwitz matrix:
  \begin{align}
    A &= \begin{bmatrix}
           -a_1 & b_{1,2} & \ldots & & b_{1,N}\\
           b_{2,1} & -a_2 & \ldots & & b_{2,N} \\
            \vdots         & & \ddots & & \\
           b_{N,1} & \ldots & \ldots & & -a_N
         \end{bmatrix}
        \label{eq:A_defn}
  \end{align}

\noindent{\bf Proof of Theorem~\ref{thm:N_feedback_same_blowup}:} Defining
$\overline V=[V_1,\ldots,V_N]^T$, we have $\dot{\overline V} \leq_e
\varphi(t)A\overline V$. Since $A$ is a strict Hurwitz matrix, a
symmetric positive-definite matrix $P$ exists such that $PA+A^TP=-Q$ with $Q$
being a symmetric positive-definite matrix. Defining $V=\overline V^TP\overline
V$, we obtain\footnote{If $P$ is a symmetric positive-definite matrix, then
  $\lambda_{max}(P)$ denotes its maximum eigenvalue and $\lambda_{min}(P)$ denotes its minimum eigenvalue.} $\dot V \leq -\varphi(t)\overline V^T Q \overline V\leq
-\varphi(t)\frac{\lambda_{min}(Q)}{\lambda_{max}(P)} V$. Hence, from
Lemma~\ref{lemma:V_converge}, $V$ is a prescribed-time exponentially convergent signal. From the
definition of $V$, we know that $V_i\leq \sqrt{\frac{V}{\lambda_{min}(P)}}, i=1,\ldots,N$.
Therefore,  $V_i,i=1,\ldots,N$ are also prescribed-time
exponentially convergent signals. Since $V_i$ is a prescribed-time exponentially
convergent Lyapunov certicate for the $i^{th}$ system, it follows that each of
the systems ($i=1,\ldots,N$) is PT-C. Therefore, the overall interconnected
system is also PT-C.
\end{thrm}

\begin{thrm}
  \label{thm:N_feedback_different_blowup}
  Consider the feedback interconnection of $N$ systems with the dynamics of the $i^{th}$ system as shown in \eqref{eq:N_feedback_interconnection_i}. If   prescribed-time exponentially convergent Lyapunov certificates $V_i(x_i,t)$
  for each of the systems ($i=1,\ldots,N$) exist satisfying
  \begin{align}
    \dot V_i &\leq \varphi_i(t) [-a_iV_i + \sum_{\substack{j=1\\j\neq i}}^{j=N} b_{i,j} V_j] \ , \ i=1,\ldots,N
               \label{eq:N_feedback_different_Vi}
  \end{align}
  with $\varphi_i,i=1,\ldots,N$ being polynomially bounded blow-up functions
  satisfying
  $\varphi_i(t) \geq
  p_{0i}\left(\frac{T}{T-t}\right)^2$ with $p_{0i},i=1,\ldots,N$ being positive constants,
  and with $a_i, i=1,\ldots,N$ and $b_{i,j}, i=1,\ldots,N, j=1,\ldots,N, j\neq
  i$, being positive constants, then the overall system formed by the feedback
  interconnection is PT-C if the matrix $A$ defined in \eqref{eq:A_defn} is such
  that $-A$ is a Lyapunov diagonally stable matrix.
  % \begin{align}
  %   A \mbox{ is strict Hurwitz} \\
  %   \lim_{t\rightarrow T}\varphi_i(t)e^{-\delta(A)\int_0^t\underline\varphi(s)ds} = 0 \ , \ i=1,\ldots,N
  % \end{align}

  \noindent{\bf Proof of Theorem~\ref{thm:N_feedback_different_blowup}:}
  Since $-A$ is a Lyapunov diagonally stable matrix, we have $\delta(A)>0$ and
  define $q=q_\delta(A)$. Analogously to the proof of
  Theorem~\ref{thm:feedback_different_blowup}, define
\begin{align}
  \overline V &= \left[
                \frac{V_1}{\varphi_1(t)} , \ldots,\frac{V_N}{\varphi_N(t)}\right]^T \ , \ V_c = [V_1,\ldots,V_N]^T.
                \label{eq:Vdefn_feedback_N}
\end{align}
Defining
$\chi=\diag(\varphi_1(t),\ldots,\varphi_N(t))$, we see that 
$V_c=\chi \overline V$.
Since $\frac{d}{dt}\varphi_i(t),i=1,\ldots,N$ are
non-negative for all time since $\varphi_i(t)$ are
monotonically increasing by the definition of a blow-up function, we see
that 
  $\dot{\overline V} \leq_e AV_c$ implying that  
  $\dot{\overline V} \leq_e A\chi \overline V$.
  Defining $V=q\overline V$ and $\underline\varphi(t)=\min\{\varphi_1(t),\ldots,\varphi_N(t)\}$, we note as in the proof of
  Theorem~\ref{thm:feedback_different_blowup} that since $q>_e 0$ and $\overline
  V\geq_e 0$, we have
$\dot V  \leq qA\chi \overline V\leq -\delta(A)
q\chi\overline V\leq -\delta(A) q\underline\varphi(t)\overline V$
implying that
$\dot V \leq -\delta(A)\underline\varphi(t)V$ since
$\delta(A)$ and $\underline\varphi(t)$ are scalars. By the Gr\"onwall-Bellman
inequality, this implies that
$V\leq e^{-\delta(A)\int_0^t\underline\varphi(s)ds}V_0$ where $V_0$ is the
value of $V$ at time $t=0$. Hence, $V$ is a
prescribed-time exponentially convergent signal. Furthermore, since $q>_e0$, we
have $V_i\leq \frac{V\varphi_i(t)}{q_i}, i=1,\ldots,N$, where $q_i$ is the
$i^{th}$ element of $q$. This yields the inequality
$V_i \leq \frac{1}{q_i} \varphi_i(t)e^{-\delta(A)\int_0^t\underline\varphi(s)ds}V_0 , i=1,\ldots,N$.
Noting that
$\varphi_i,i=1,\ldots,N$ are polynomially bounded blow-up functions and
$\underline\varphi(s)ds\geq \underline
p_0\left(\frac{T}{T-t}\right)^2$ where $\underline
p_0=\min\{p_{01},\ldots,p_{0N}\}$,  we see that $V_i$ is upper bounded by a product of
a polynomially bounded blow-up function and a prescribed-time exponentially
convergent signal. By Lemma~\ref{lemma:polynomial_exp_convergence2}, this
implies that $V_i, i=1,\ldots,N$ are prescribed-time exponentially convergent signals.
Since $V_i$ is a
prescribed-time exponentially convergent Lyapunov certificate for the $i^{th}$
system, it follows that all the systems ($i=1,\ldots,N$) and the overall interconnected system are PT-C.
\end{thrm}

% \begin{thrm}
%   \label{thm:N_feedback_different_blowup_poly}
%   Consider the feedback interconnection of $N$ systems with the dynamics of the $i^{th}$ system as shown in \eqref{eq:N_feedback_interconnection_i}. If positive-definite functions $V_i(x_i,t), i=1,\ldots,N$, exist satisfying \eqref{eq:N_feedback_different_Vi} 
%   with $\varphi_i,i=1,\ldots,N$ being blow-up functions and $a_i, i=1,\ldots,N$ and $b_{i,j}, i=1,\ldots,N, j=1,\ldots,N, j\neq i$, being positive constants, then the overall system formed by the feedback interconnection is PT-C if 
%     the matrix $A$ defined in \eqref{eq:A_defn} is strict Hurwitz.
% \end{thrm}
% \noindent{\bf Proof of Theorem~\ref{thm:N_feedback_different_blowup_poly}:}

\vspace*{-0.2in}
\section{Illustrative Examples}
\label{sec:examples}
To illustrate the results in Section~\ref{sec:interconnections}, we consider two
examples in this section, the first example being a simple cascade
interconnection of two given scalar systems and the second example being a
decentralized control design problem for a feedback interconnection of two
systems.

\vspace*{-0.1in}
\subsection{Example 1: Cascade Interconnection}
\begin{figure}[!h]
  \centering
  \ifdef{\TwoColumn}{\includegraphics[width=2.1in,clip=true,trim=0in 0.07in 0in 0.35in]{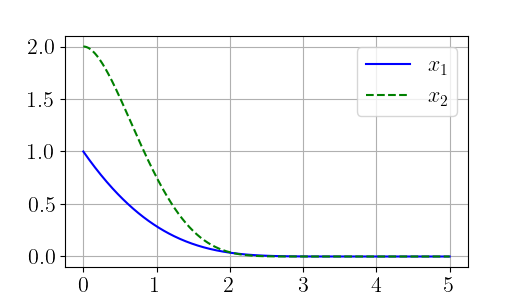}}{\includegraphics[width=2.9in,clip=true,trim=0in 0.07in 0in 0.35in]{sim1_cascade.png}}
  \caption{Simulation results for Example 1 (cascade interconnection).}
  \label{fig:sim_cascade}
  \vspace*{-0.2in} 
\end{figure}

Consider the cascade interconnection of the following two scalar systems (the state of
the first system being $x_1$ and the state of the second system being $x_2$):
\begin{align}
  \dot x_1 &= -\varphi_1(t)x_1 \ \ \ ; \ \ \ 
  \dot x_2 = -\varphi_2(t)x_2 + (1+\varphi_2^2(t))x_1^3
             \label{eq:cascade_example_sys12}
\end{align}
where $\varphi_1(t)$ and $\varphi_2(t)$ are blow-up functions given as
$\varphi_1(t)=\left(\frac{T}{T-t}\right)^2$ and
$\varphi_2(t)=\left(\frac{T}{T-t}\right)^3$. Defining the Lyapunov functions
$V_1=\frac{1}{2}x_1^2$ and $V_2=\frac{1}{2}x_2^2$, we note that
\begin{align}
  \dot V_1 &= -\varphi_1(t)x_1^2 = -2\varphi_1(t)V_1 \\
  \dot V_2 &= -\varphi_2(t)x_2^2 + (1+\varphi_2^2(t))x_1^3x_2 \nonumber\\
           %&\leq -\frac{1}{2}\varphi_2(t)x_2^2 + \frac{(1+\varphi_2^2(t))^2}{2\varphi_2(t)}x_1^6\nonumber\\
           &\leq 
             -\varphi_2(t)V_2 + \frac{4(1+\varphi_2^2(t))^2}{\varphi_2(t)}V_1^3.
\end{align}
The conditions of Theorem~\ref{thm:cascade_different_blowup} are seen to be
satisfied with $a=1$, $p(V_1)=V_1^3$, and
$\varphi_3(t)=\frac{4(1+\varphi_2^2(t))^2}{\varphi_2^2(t)}$. Hence, by Theorem~\ref{thm:cascade_different_blowup}, the cascade interconnection shown in equations
\eqref{eq:cascade_example_sys12} is PT-C.
This is verified in the simulation in Figure~\ref{fig:sim_cascade}. The initial
conditions for the simulation are picked to be $x_1=1$ and $x_2=2$.
The prescribed time $T$ is set to be 5~s. To avoid numerical issues in
simulations, we set an effective terminal time $\overline T$ to be a slightly
larger constant than $T$ as per the standard procedure in prior works \cite{prescribed_time3,KKK20a}: $\overline T=5.05$~s.
It is seen
in Figure~\ref{fig:sim_cascade} that the states of both systems converge to 0 as
$t\rightarrow T$.

\vspace*{-0.1in}
\subsection{Example 2: Feedback Interconnection}
Consider the feedback interconnection of two systems as shown below:
\begin{itemize}
\item System 1 (with state $[x_{11},x_{12}]^T$ and input $u_1$):
  \begin{align}
    \dot x_{11} &= x_{12} + x_{21}^3\ \ \ ; \ \ \ 
    \dot x_{12} = u_1.
  \end{align}
\item System 2 (with state $[x_{21},x_{22}]^T$ and input $u_2$):
  \begin{align}
    \dot x_{21} &= x_{22} + \sin(x_{21})x_{11} \ \ \ ; \ \ \ 
    \dot x_{22} = u_2.
  \end{align}
\end{itemize}
Consider that we wish to design a decentralized controller
(i.e., $u_1$ having state dependence only on $[x_{11},x_{12}]^T$; $u_2$ having
state dependence only on $[x_{21},x_{22}]^T$) to make the overall interconnected
system PT-C. We design prescribed-time controllers for the two systems
below using two different control design methods and then
show as an application of Theorem~\ref{thm:feedback_different_blowup} that the
interconnected system is PT-C.

For the first system, we use a backstepping approach. Define the virtual
control law $x_{12}^*=-k_{11}\varphi_1(t)x_{11}-k_{12}\varphi_1(t)x_{11}^9 $
with $k_{11}$ and $k_{12}$ being positive constants
and with $\varphi_1(t)$ being a polynomially bounded blow-up function to
be picked below.
The choice of the term involving $x_{11}^9$ is motivated by the fact noted below
that the upper bounding of the expression obtained in the Lyapunov analysis
from the interconnection term $x_{21}^3$ will generate a term involving $x_{11}^{12}$.
Then, with
$z_{12}=x_{12}-x_{12}^*$, we define the Lyapunov function
$V_1=\frac{1}{4}x_{11}^4+\frac{c_1}{2}z_{12}^2$ with $c_1$ being a positive
constant.
Then, we have
\begin{align}
  \dot V_1 &= x_{11}^3(x_{12}+x_{21}^3)+c_1z_{12}(u_1 +\dot\varphi_1(t)(k_{11}x_{11}+k_{12}x_{11}^9)
             \nonumber\\
         &\quad    +\varphi_1(t)(k_{11}+9k_{12}x_{11}^8)(x_{12}+x_{21}^3)).
             \label{eq:V1dot_ex_feedback}
\end{align}
Using the inequality $ab\leq \frac{1}{p}|a|^p + \frac{1}{q}|b|^q$ which holds for
any real numbers $a$ and $b$ and any positive constants $p$ and $q$ with $\frac{1}{p}+\frac{1}{q}=1$, we can write
the following inequalities to upper bound the various ``interconnection'' terms
(i.e., terms involving state variables of both System 1 and System 2)
appearing in the right hand side of \eqref{eq:V1dot_ex_feedback}:
\begin{align}
  x_{11}^3x_{21}^3 &\leq \frac{3}{4(4\varphi_1(t)k_{12})^{\frac{1}{3}}}x_{21}^4+k_{12}\varphi_1(t)x_{11}^{12}
                     \label{eq:ineq1_ex_feedback}
\end{align}
\begin{align}
  c_1z_{12}\varphi_1(t)(k_{11}+9k_{12}x_{11}^8)x_{21}^3 &\leq \frac{3}{4(4\varphi_1(t)k_{12})^{\frac{1}{3}}}x_{21}^4
                                                          \nonumber\\
                                                      &\!\!\!\!\!\!\!\!\!\!\!\!\!\!\!
                                                          +k_{12}c_1^4z_{12}^4\varphi_1^5(t)(k_{11}+9k_{12}x_{11}^8)^4.
                                                          \label{eq:ineq2_ex_feedback}
\end{align}
Designing the control input $u_1$ as
\begin{align}
  u_1 &= -k_{13}z_{12}\varphi_1(t)-\frac{1}{c_1}x_{11}^3 -\dot\varphi_1(t)(k_{11}x_{11}+k_{12}x_{11}^9)
        \nonumber\\
      &\quad
        -\varphi_1(t)(k_{11}+9k_{12}x_{11}^8)x_{12}
        \nonumber\\
  &\quad
        -k_{12}c_1^3z_{12}^3\varphi_1^5(t)(k_{11}+9k_{12}x_{11}^8)^4
        \label{eq:u1_ex_feedback}
\end{align}
with $k_{13}$ being a positive constant, 
and using the inequalities \eqref{eq:ineq1_ex_feedback} and
\eqref{eq:ineq2_ex_feedback}, \eqref{eq:V1dot_ex_feedback} reduces to
\begin{align}
  \dot V_1 &\leq -k_{11}\varphi_1(t)x_{11}^4 -k_{13} \varphi_1(t) c_1z_{12}^2 + \frac{3}{2(4\varphi_1(t)k_{12})^{\frac{1}{3}}}x_{21}^4.
             \label{eq:V1dot2_ex_feedback}
\end{align}

To design the controller for the second system (i.e., the system with state
variables $[x_{21},x_{22}]^T$), we apply a high-gain scaling analogous to
\cite{KKK19a,KKK20a}, but with a polynomially bounded blow-up function $\varphi_2(t)$ in
place of a dynamic high-gain scaling parameter $r$:
\begin{align}
  \overline x_{21} &= x_{21} \ ; \ \overline x_{22} = \frac{x_{22}}{\varphi_2(t)}.
\end{align}
In terms of these scaled state variables, we have the dynamics
\begin{align}
  \dot{\overline x}_{21} &= \varphi_2(t)\overline x_{22} + \sin(\overline x_{21}) x_{11}
                           \nonumber\\
  \dot{\overline x}_{22} &= -\frac{\dot\varphi_2(t)}{\varphi_2(t)}\overline x_{22}
                           + \frac{1}{\varphi_2(t)} u_2.
\end{align}
Defining the Lyapunov function of form
\begin{align}
  V_2 &= \frac{1}{4}\overline x_{21}^4 +
        \frac{c_2}{2}(\overline x_{22}+k_{21}\overline x_{21})^2
\end{align}
with $k_{21}$ being a
positive constant, we obtain
\begin{align}
  \dot V_2 &= \overline x_{21}^3(\varphi_2(t)z_{22} -\varphi_2(t)k_{21}\overline x_{21} + \sin(\overline x_{21}) x_{11})
\nonumber\\
&\quad             + c_2z_{22}\Big(-\frac{\dot\varphi_2(t)}{\varphi_2(t)}\overline x_{22}
  + \frac{1}{\varphi_2(t)} u_2
  \nonumber\\
  &\quad
  + k_{21}(\varphi_2(t)\overline x_{22} + \sin(\overline x_{21}) x_{11})\Big)
             \label{eq:V2dot_ex_feedback}
\end{align}
where $z_{22}=\overline x_{22}+k_{21}\overline x_{21}$. We can write the
following inequalities to upper bound the various interconnection terms appearing in the right
hand side of \eqref{eq:V2dot_ex_feedback}:
\begin{align}
  \overline x_{21}^3 \sin(\overline x_{21}) x_{11} &\leq \frac{1}{4}\varphi_2(t)k_{21}\overline x_{21}^4 + \frac{27}{4(\varphi_2(t)k_{21})^3}x_{11}^4
                                                     \label{eq:ineq3_ex_feedback}
                                                     \\
  c_2 z_{22}k_{21}\sin(\overline x_{21}) x_{11} &\leq \frac{1}{2}c_2^2 k_{21}^2 z_{22}^2  + \frac{1}{4}\varphi_2(t)k_{21}\overline x_{21}^4 +
                                                   \nonumber\\
                                                    &\quad
                                                      \frac{1}{4\varphi_2(t)k_{21}}x_{11}^4.
                                                      \label{eq:ineq4_ex_feedback}
\end{align}
Designing the control input $u_2$ as
\begin{align}
  u_2 &= -\dot\varphi_2(t)k_{21}\overline x_{21}+\varphi_2(t)\Big(
        -k_{21}\varphi_2(t)\overline x_{22} - \frac{1}{c_2}\varphi_2(t)\overline x_{21}^3
\nonumber\\
&\quad        -\frac{1}{2}c_2k_{21}^2z_{22} - \varphi_2(t) k_{22} z_{22}
  \Big)
  \label{eq:u2_ex_feedback}
\end{align}
with $k_{22}$ being a positive constant, using the inequalities \eqref{eq:ineq3_ex_feedback} and
\eqref{eq:ineq4_ex_feedback}, and noting that
$\dot\varphi_2(t)$ is non-negative for any blow-up function $\varphi_2$,
\eqref{eq:V2dot_ex_feedback} reduces to
\begin{align}
  \dot V_2&\leq -\frac{1}{2}\varphi_2(t)k_{21}\overline x_{21}^4 - k_{22}\varphi_2(t)c_2 z_{22}^2
            \nonumber\\
  &\quad
    + \Big(
\frac{27}{4(\varphi_2(t)k_{21})^3}+\frac{1}{4\varphi_2(t)k_{21}}
    \Big)x_{11}^4.
    \label{eq:V2dot2_ex_feedback}
\end{align}
Denoting $\varphi_{10}=\varphi_1(0)$ and $\varphi_{20}=\varphi_2(0)$, we see
from \eqref{eq:V1dot2_ex_feedback} and \eqref{eq:V2dot2_ex_feedback} that
inequalities of the form \eqref{eq:feedback_different_V1} and
\eqref{eq:feedback_different_V2} are satisfied
with
$a_1=\min(4k_{11},2k_{13})$,
$a_2=\min(2k_{21},2k_{22})$,
$b_1=\frac{6}{\varphi_{10}(4\varphi_{10}k_{12})^{\frac{1}{3}}}$,
and $b_2=\frac{27}{\varphi_{20}(\varphi_{20}k_{21})^{3}}+\frac{1}{\varphi_{20}^2k_{21}}$.
From Theorem~\ref{thm:feedback_different_blowup}, it follows that the feedback
interconnection of the two systems is PT-C if the controller parameters are
picked such that $a_1a_2>b_1b_2$. Note that the control designs for $u_1$ and
$u_2$ in \eqref{eq:u1_ex_feedback} and \eqref{eq:u2_ex_feedback} are based on
two different methods and are
decentralized in the sense that each control law depends only on the
state variables of the corresponding system (i.e., $u_1$ is a function of
$[x_{11},x_{12}]^T$; $u_2$ is a function of $[x_{21},x_{22}]^T$).

For simulation studies, we pick the blow-up functions for the two systems to be
$\varphi_1(t)=6+\left(\frac{T}{T-t}\right)^2$ and
$\varphi_2(t)=6+\left(\frac{T}{T-t}\right)^3$, respectively. The motivation for
including positive constants in the choices of $\varphi_1(t)$ and $\varphi_2(t)$
is to increase the values of $\varphi_{10}=\varphi_1(0)$ and $\varphi_{20}=\varphi_2(0)$ taking into
account the appearance of these values in the expressions for $b_1$ and $b_2$ so
as to make $b_1$ and $b_2$ smaller making it more likely that the condition
$a_1a_2>b_1b_2$ is satisfied.
The controller parameters are picked to be
$k_{11}=0.25$,
$k_{12}=0.1$,
$k_{13}=0.5$,
$c_1=0.05$,
$k_{21}=0.25$,
$k_{22}=0.5$, and $c_2=0.02$. With these choices of the controller parameters,
we see that the parameters $a_1$, $a_2$, $b_1$, and $b_2$ in
\eqref{eq:feedback_different_V1} and \eqref{eq:feedback_different_V2} are
$a_1=1$, $a_2=0.5$, $b_1=0.608$, and $b_2=0.801$. With these parameters, we have
$a_1a_2=0.5$ and $b_1b_2=0.487$. Hence, the condition
$a_1a_2>b_1b_2$ is satisfied implying that the feedback interconnection of Systems 1 and 2 is
PT-C by Theorem~\ref{thm:feedback_different_blowup}.  

\begin{figure}[!h]
  \centering
  \ifdef{\TwoColumn}{\includegraphics[width=2.8in,clip=true,trim=0in 0.25in 0in 0.8in]{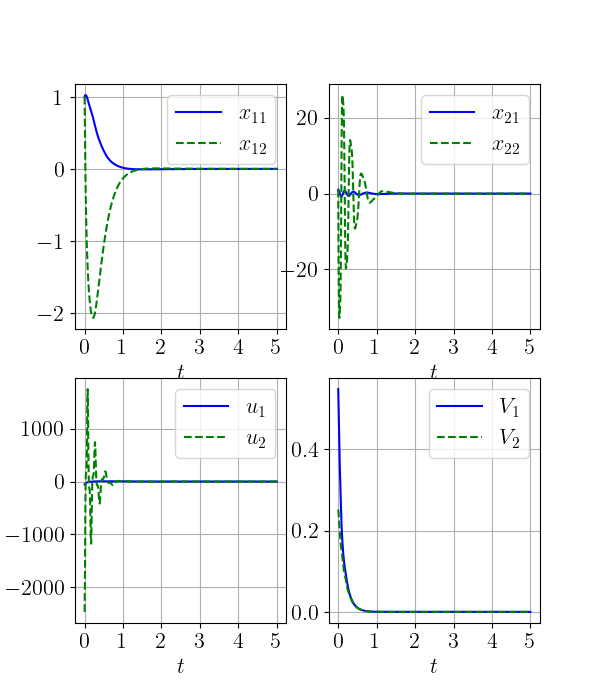}}{\includegraphics[width=3.5in,clip=true,trim=0in 0.25in 0in 0.8in]{sim1_feedback.png}}
  \caption{Simulation results for Example 2 (feedback interconnection).}
  \label{fig:sim_feedback}
  \vspace*{-0.2in} 
\end{figure}

Simulation results for the closed-loop interconnected system using the 
control designs described above are shown in
Figure~\ref{fig:sim_feedback}.
The initial conditions for the simulation are specified as
$x_{11}=x_{12}=x_{21}=x_{22}=1$.  
As in Example 1, the prescribed time $T$ is set
to be 5~s and the effective terminal time $\overline T$ to avoid numerical
issues in simulations is set to be $\overline T=5.05$.
It is seen that all the
state variables as well as the control inputs converge to 0 as $t\rightarrow T$.
Note however that the
condition $a_1a_2>b_1b_2$ is a sufficient but not necessary condition and
furthermore, there is some intrinsic conservativeness in the upper bounds 
\eqref{eq:ineq1_ex_feedback}, \eqref{eq:ineq2_ex_feedback},
\eqref{eq:ineq3_ex_feedback}, and \eqref{eq:ineq4_ex_feedback}. Hence, the
controller parameters could be tuned to reduce overshoots (e.g., picking 
smaller controller gains such as $k_{11}=k_{21}=0.1$ and removing the additive
positive constants in the choices of $\varphi_1(t)$ and $\varphi_2(t)$ results
in smoother closed-loop trajectories with greatly reduced values of the control
signals). Simulation plots with these modified controller parameters are omitted
for brevity. 

% k11 = 0.1
% k12 = 0.1
% k13 = 0.5
% c1 = 0.05
% k21 = 0.1
% k22 = 0.5
% c2 = 0.02

In the control designs above, note that $u_1$ and $u_2$ involve the time
derivatives $\dot\varphi_1(t)$ and $\dot\varphi_2(t)$, respectively. The
convergence of $u_1$ and $u_2$ to zero inspite of this dependence on the time
derivatives of the blow-up functions is due to the fact that these time
derivatives $\dot\varphi_1(t)$ and
$\dot\varphi_2(t)$ are also polynomially bounded blow-up functions. Hence,
the terms in $u_1$ and $u_2$ involving these time derivatives are products of
polynomially bounded blow-up functions and prescribed-time exponentially
convergent signals and it is known from
Lemma~\ref{lemma:polynomial_exp_convergence2} that such products are also prescribed-time exponentially
convergent (to zero). To see that the time derivatives $\dot\varphi_1(t)$ and
$\dot\varphi_2(t)$ are also polynomially bounded blow-up functions, note that for any blow-up function $\varphi(t)$ which is a polynomial in
$\frac{T}{T-t}$, its time derivative $\dot\varphi(t)$ is also a polynomial in
$\frac{T}{T-t}$ since
$\frac{d}{dt}\left(\frac{T}{T-t}\right)^a=\frac{a}{T-t}\left(\frac{T}{T-t}\right)^a$
for any $a\geq 1$. Hence, $\dot\varphi(t)$ is also a polynomially bounded blow-up function.

 \section{Conclusion}
 \label{sec:conclusion}
 Cascade and feedback interconnections of prescribed-time ISS systems were considered and sufficient conditions were
 developed under which such interconnections retain prescribed-time
 convergence properties. 
 Interconnections both of two systems and of an arbitrary number of systems were
 considered.
 Central tools in the analysis were the newly introduced notions of polynomially
 bounded blow-up functions and prescribed-time exponential convergence and the
 related concept of prescribed-time exponentially convergent Lyapunov certificates. A
 detailed analysis of these new notions was presented and their links to
 prescribed-time convergence properties of interconnected systems were explored.
 Simulation studies were performed for example systems in cascade and feedback
 interconnection structures. 
 In the feedback interconnection example, a scenario where controllers are
 designed separately for the two systems in the interconnection (indeed, using
 two different control design methods) was considered and it was shown as an
 application of the developed theoretical results that the controllers can be
 put together to achieve prescribed-time stabilization of the interconnected
 system. The newly introduced notions noted above provide powerful tools for
 analyzing prescribed-time properties of both single systems and interconnected
 systems. Future work will address further study of the application of these
 notions in the context of more general interconnection structures (e.g.,
 general nonlinear and time-varying dependencies in Lyapunov inequalities of each system on Lyapunov certificates of other systems) involving wider classes of
 nonlinear systems and control designs
 including adaptive and output-feedback controllers.

 %It was shown as an application of
 %the developed theoretical results on interconnections of prescribed-time ISS
 %systems that dynamic controllers designed separately for each of the two
 %coupled nonlinear systems can be put together to achieve prescribed-time
 %stabilization of the interconnected system.

%%\vphantom{\cite{KKK95,Isi99,Kha01,KA01,Tee92,MP96,SJK97b,PW96,ZG06,RGLS08,Jan05,MB06,MML07,KKC02,KK04e,KK04f,KS87,IR03,Lin09,Ito06,Pra03,KK02d,KAP06,KK07f,KK07a,KK07h,KK10b,KK10a,KK11a,KK13a,SOK84,HLS07,Yur08,CA09}}

%% \renewcommand{\baselinestretch}{1.0}
 \bstctlcite{IEEEexample:BSTcontrol}
 \bibliographystyle{IEEEtran}
 \bibliography{refs}

\end{document}